\documentclass[final]{siamart220329}  

\usepackage{siampaper}

\usepackage{tikz}
\usepackage{pgfplots}
\usepackage{tikzscale}
\usetikzlibrary{external}
\def\figurespath{figures}
\graphicspath{{\figurespath/}}
\tikzset{external/system call={lualatex \tikzexternalcheckshellescape
      -halt-on-error
      -interaction=batchmode
      -jobname "\image" "\texsource"}}

\newcommand*{\includetikzgraphics}[2][]{%
      \includegraphics[#1]{#2}
}

\makeatletter
\renewcommand{\todo}[2][]{\tikzexternaldisable\@todo[#1]{#2}\tikzexternalenable}
\newcommand*{\defeq}{\mathrel{\vcenter{\baselineskip0.5ex \lineskiplimit0pt
                     \hbox{\scriptsize.}\hbox{\scriptsize.}}}%
                     =}
\makeatother

\newtheorem{problemassumption}{Problem Assumption}[section]
\newtheorem{modelassumption}{Model Assumption}[section]

\crefname{problemassumption}{Problem Assumption}{Problem Assumptions}
\Crefname{problemassumption}{Problem Assumption}{Problem Assumptions}
\crefname{modelassumption}{Model Assumption}{Model Assumptions}
\Crefname{modelassumption}{Model Assumption}{Model Assumptions}
\crefname{stepassumption}{Step Assumption}{Step Assumptions}
\Crefname{stepassumption}{Step Assumption}{Step Assumptions}

\newcommand{\dom}{\mathop{\mathrm{dom}}}
\newcommand{\varphicp}{\varphi_{\textup{cp}}}
\newcommand{\xicp}{\xi_{\textup{cp}}}
\newcommand{\xikcp}{\xi_{k, \textup{cp}}}
\newcommand{\scp}{s_{\textup{cp}}}
\newcommand{\skcp}{s_{k, \textup{cp}}}
\newcommand{\Kf}{\mathcal{K}_f}

\newcommand{\cahiernumber}{64}  

\def\papertitle{A Proximal Modified Quasi-Newton Method for Nonsmooth Regularized Optimization}
\def\authorone{YOUSSEF DIOUANE}
\def\authortwo{Mohamed L. Habiboullah}
\def\authorthree{Dominique Orban}

\hypersetup{
  pdftitle={\papertitle},
  pdfauthor={\authorone{}, \authortwo{} and \authorthree{}},
  pdfsubject={nonsmooth regularized optimization},
  pdfkeywords={Nonsmooth optimization, Nonconvex optimization, Regularized optimization, Composite optimization, Modified quasi-Newton method, Proximal quasi-Newton method, Proximal gradient method},
}

\title{\papertitle}

\author{
  \authorone\thanks{%
    GERAD and Department of Mathematics and Industrial Engineering, Polytechnique Montr\'eal.
    E-mail: \mailto{youssef.diouane@polymtl.ca}.
    Research partially supported by an NSERC Discovery Grant.}
  \and
  \authortwo\thanks{%
    GERAD and Department of Mathematics and Industrial Engineering, Polytechnique Montr\'eal.
    E-mail: \mailto{mohamed-laghdaf-2.habiboullah@polymtl.ca}.
    Research partially supported by an FRQNT scholarship.
  }
  \and
  \authorthree\thanks{%
    GERAD and Department of Mathematics and Industrial Engineering, Polytechnique Montr\'eal.
    E-mail: \mailto{dominique.orban@gerad.ca}.
    Research partially supported by an NSERC Discovery Grant.
  }
}
\date{\today}

\begin{document}

\maketitle

\thispagestyle{firstpage}
\pagestyle{myheadings}

\begin{abstract}
  We develop R2N, a modified quasi-Newton method for minimizing the sum of a \(\mathcal{C}^1\) function \(f\) and a lower semi-continuous prox-bounded \(h\).
  Both \(f\) and \(h\) may be nonconvex.
  At each iteration, our method computes a step by minimizing the sum of a quadratic model of \(f\), a model of \(h\), and an adaptive quadratic regularization term.
  A step may be computed by way of a variant of the proximal-gradient method.
  An advantage of R2N over competing trust-region methods is that proximal operators do not involve an extra trust-region indicator.
  We also develop the variant R2DH, in which the model Hessian is diagonal, which allows us to compute a step without relying on a subproblem solver when \(h\) is separable.
  R2DH can be used as standalone solver, but also as subproblem solver inside R2N\@.
  We describe non-monotone variants of both R2N and R2DH\@.
  Global convergence of a first-order stationarity measure to zero holds without relying on local Lipschitz continuity of \(\nabla f\), while allowing model Hessians to grow unbounded, an assumption particularly relevant to quasi-Newton models.
  Under Lipschitz-continuity of \(\nabla f\), we establish a tight worst-case evaluation complexity bound of \(O(1 / \epsilon^{2/(1 - p)})\) to bring said measure below \(\epsilon > 0\), where \(0 \leq p < 1\) controls the growth of model Hessians.
  Specifically, the latter must not diverge faster than \(|\mathcal{S}_k|^p\), where \(\mathcal{S}_k\) is the set of successful iterations up to iteration \(k\).
  When \(p = 1\), we establish the tight exponential complexity bound \(O(\exp(c \epsilon^{-2}))\) where \(c > 0\) is a constant.
  We describe our Julia implementation and report numerical experience on a classic basis-pursuit problem, an image denoising problem, a minimum-rank matrix completion problem, a nonlinear support vector machine and an inverse nonlinear problem.
\end{abstract}



\section{Introduction}%
\label{sec:introduction}

We consider problems of the form
\begin{equation}%
  \label{eq:nlp}
  \minimize{x \in \R^n} \ f(x) + h(x),
\end{equation}
where \(f: \R^n \to \R\) is \(\mathcal{C}^1\) on \(\R^n\), and \(h: \R^n \to \R \cup \{+\infty\}\) is lower semi-continuous (lsc).
Both \(f\) and \(h\) may be nonconvex.
Our motivation is to develop a modified Newton variant of the trust-region algorithms of \citet{aravkin-baraldi-orban-2022} and \citet{leconte-orban-2025} because the proximal operators used in the subproblems should be easier to derive.
For instance, when \(h\) is the rank function, the proximal operator with a trust-region indicator is not known analytically at this time.

We introduce method R2N, at each iteration of which the sum of a quadratic model of \(f\), a model of \(h\), and an adaptive quadratic regularization term, is approximately minimized.
Both models may be nonconvex.
The Hessian of the quadratic model of \(f\) may be that of \(f\) if it exists, or an approximation such as those derived from quasi-Newton updates.
We establish global convergence of R2N under the assumption that the models of \(h\) are prox-bounded and approximate \(h(x + s)\) as \(o(\|s\|)\)---an assumption that covers composite terms with H\"older Jacobian, see \Cref{asm:model-h-adequacy} for details.
No assumption on local Lipschitz continuity of \(\nabla f\) is required, nor is boundedness of the model Hessians, provided they do not diverge too fast.
Specifically, if \(B_k\) is the model Hessian at iteration \(k\), we require that the series with general term \(1 / (1 + \max_{0 \leq j \leq k} \|B_j\|)\) diverge---an assumption similar to that used in trust-region methods \citep[\S\(8.4\)]{conn-gould-toint-2000}.
Our assumptions are significantly weaker than assumptions commonly found in the analysis of competing methods, and, consequently, the applicability of R2N is significantly more general---see the related research section below for details.

R2N specializes to method R2DH when \(B_k\) is diagonal, as did the solver of \citet{leconte-orban-2025}.
For a number of choices of separable \(h\) that are relevant in applications, steps can be computed explicitly without resort to an iterative subproblem solver.
R2DH can be used as standalone solver or as subproblem solver inside R2N.

We also develop complexity results inspired from those of \citet{leconte-orban-2023-2} and \citet{diouane-habiboullah-orban-2024}, that account for potentially unbounded model Hessians.
Specifically, we require that either \(\|B_k\| = O(|\mathcal{S}_k|^p)\) for some \(0 \leq p \leq 1\), where \(\mathcal{S}_k\) is the set of successful iterations up to iteration \(k\).
When \(0 \leq p < 1\), we establish a tight \(O(\epsilon^{-2/(1 - p)})\) complexity, and when \(p = 1\), we establish a tight exponential complexity, i.e., a bound in \(O(\exp(c \epsilon^{-2}))\) where \(c > 0\) is a constant.
Though the latter bound is tight, it is not known if it is attained for a quasi-Newton update.

We provide efficient implementations of R2N and R2DH\@.
The latter can use one of several diagonal quasi-Newton updates.
Both have non-monotone variants that preserve their convergence and complexity properties.
Our open-source Julia implementations are available from \citep{baraldi-leconte-orban-regularized-optimization-2024}.
In \Cref{sec:numerical}, we illustrate the performance of R2N and R2DH on challenging problems, including minimum-rank problems for which the trust-region methods of \citep{aravkin-baraldi-orban-2022,aravkin-baraldi-orban-2024} are impractical.

\subsection*{Contributions and related research}

The proximal-gradient method \citep{fukushima-mine-1981,lions-mercier-1979} is the prototypical first-order method for~\eqref{eq:nlp}.
Vast amounts of literature consider variants but restrict \(f\) and/or \(h\) to be convex, impose that \(f\) have (locally) Lipschitz-continuous gradient, or that \(h\) be Lipschitz continuous.
For instance, \citep{lee-sun-saunders-2014} develop a proximal Newton method that requires \(f\) and \(h\) convex, a positive semi-definite Hessian, and solve the subproblem via the proximal-gradient method.
\citet{cartis-gould-toint-2011} require \(h\) to be globally Lipschitz continuous.
\citet{kanzow-lechner-2024,liu-pan-wu-yang-2024} develop an approach closely related to ours, but for convex \(h\).
Others dispense with convexity but require coercivity of \(f + h\) \citep{li-lin-2015}.

We are aware of several references that allow both \(f\) and \(h\) to be nonconvex.
For instance, \citet{bolte-sabach-teboulle-2014} propose PALM, an alternating first-order method for problems with partitioned variables.
They assume that \(f\), which acts on both sets of variables \(x\) and \(y\), has a gradient that is Lipschitz continuous with respect to \(x\) and \(y\) separately, and is Lipschitz continuous in \((x, y)\) on a bounded set.
They also assume that the sequence generated by their method is bounded, which is a strong requirement.
Under those assumptions, they prove that every accumulation point of the iterates is a first-order stationary point.
Moreover, if the K\L{} property holds, then the entire sequence converges to a first-order stationary point.
No second-order information is used in their method.

\citet{bot-csetnek-laszlo-2016} study a proximal-gradient algorithm with momentum for solving~\eqref{eq:nlp}.
They assume that \(\nabla f\) is Lipschitz continuous, that \(h\) is bounded below, and that \(f + h\) is coercive.
Under those conditions, they show that every accumulation point of the iterates is a first-order stationary point of~\eqref{eq:nlp}.
If, in addition, the function \(H(x, y) \;=\; f(x) + h(x) + M \|x - y\|^2\), where \(M > 0\), is a K\L{} function, then the entire sequence converges to a first-order stationary point.
Their method does not use second-order information, and the coercivity assumption can be restrictive.

\citet{themelis-stella-patrinos-2017} propose ZeroFPR, a non-monotone line-search proximal quasi-Newton method based on the forward-backward envelope.
They assume that \(\nabla f\) is globally Lipschitz continuous—although they note that local Lipschitz continuity suffices if the domain of \(h\) is bounded and the search directions remain bounded—and that \(h\) is prox-bounded.
Under these assumptions, any accumulation point is a stationary point of~\eqref{eq:nlp}.
Furthermore, if the iterates remain bounded, the forward-backward envelope satisfies the K\L{} property, and \(f\) is twice continuously differentiable, then the entire sequence converges to a first-order stationary point.
With a suitable desingularization function, they also obtain R-linear convergence of \(\{x_k\}_{k\in\mathbb{N}}\).
Although ZeroFPR allows quasi-Newton approximations, the model Hessians are assumed to be uniformly bounded—a condition that may be difficult to guarantee, as we will discuss later.
It furthermore requires an estimate of the Lipschitz constant of \(\nabla f\), and a preliminary loop to compute it.

\citet{stella-themelis-sopasakis-patrinos-2017} propose PANOC, a line-search limited-memory BFGS method, in the context of optimal control problems.
They assume that \(\nabla f\) is Lipschitz continuous, although local Lipschitz continuity suffices if the domain of \(h\) is bounded and the search directions remain bounded using the similar arguments as in \citep{themelis-stella-patrinos-2017}.
However, if the Lipschitz constant of \(\nabla f\) is unknown, it must be estimated in a preliminary loop.
In addition, they require \(h\) to be bounded below, which is a strong assumption.
Under these conditions, any accumulation point of the iterates is stationary for~\eqref{eq:nlp}.
Moreover, if the iterates converge to a strong local minimum of \(f + h\), the forward-backward envelope is twice continuously differentiable, the proximal operator of \(h\) is strictly differentiable, and the model Hessians \(B_k\) satisfy the Denis–Moré condition, then the entire sequence converges at a superlinear rate.

One may observe that in the methods of \citep{bolte-sabach-teboulle-2014,bot-csetnek-laszlo-2016,themelis-stella-patrinos-2017}, the K\L{} property together with the Lipschitz continuity of \(\nabla f\) is used to establish convergence of the entire sequence of iterates to a first-order stationary point.
Moreover, with a suitable desingularization function, R-linear convergence of the iterates is shown in \citep{bot-csetnek-laszlo-2016,themelis-stella-patrinos-2017}.
By contrast, our work does not rely on such strong assumptions.
We avoid both the K\L{} property and the Lipschitz continuity of \(\nabla f\), and aim instead for generality.

More recently, some works have derived convergence analyses without assuming Lipschitz continuity of \(\nabla f\) or the K\L{} property.
For example, \citet{kanzow-mehlitz-2022} analyze monotone and non-monotone first-order proximal-gradient methods under the assumptions that \(f\) is continuously differentiable and \(h\) is lower semicontinuous and bounded below by an affine function.
\citet{de-marchi-2023} extends their setting to allow \(h\) to be prox-bounded.
They establish that for any convergent subsequence \(\{x_k\}_{k\in\mathcal{K}}\to x^*\), a stationarity measure converges to zero along \(\mathcal{K}\).
Additional strong assumptions—such as local Lipschitz continuity of \(\nabla f\) or continuity of \(h\)—then ensure that each accumulation point is a first-order stationary point.
Theirs are first-order methods and do not incorporate second-order information.
Furthermore, their convergence results do not cover cases where the iterates may be unbounded.

In this work, we establish the global convergence of a first-order stationarity measure to zero under the minimal assumptions that \(f\) is continuously differentiable and \(h\) is lower semicontinuous and prox-bounded.
We further guarantee that, for any \(\epsilon > 0\), our stationarity measure falls below \(\epsilon\) in a finite number of iterations, even if the sequence of iterates is unbounded.

Our work follows the scheme laid out by \citet{aravkin-baraldi-orban-2022}; a trust-region framework applicable to nonconvex \(f\) and/or \(h\), and that does not require coercivity or K\L\, assumptions.
However, their analysis relies on the Lipschitz continuity of \(\nabla f\) in a neighborhood of the iterates.
They also describe a method named R2 that amounts to a proximal-gradient method with adaptive step size, and that may be viewed as R2N where \(B_k\) is set to zero at each iteration, effectively reducing to a first-order method.
\citet{aravkin-baraldi-orban-2024} specialize their trust-region method to problems where \(f\) has a least-squares structure, and develop a Levenberg-Marquardt variant named LM that may also be viewed as a special case of R2N for least-squares \(f\).
If \(J_k\) is the least-squares residual's Jacobian at \(x_k\), their model of \(f\) uses \(B_k = J_k^T J_k\).
\citet{leconte-orban-2025} devise variants of the trust-region method of \citep{aravkin-baraldi-orban-2022} for separable \(h\) in which the model Hessian is a diagonal quasi-Newton approximation.
They also devise non-monotone schemes that are shown to significantly improve performance in certain cases.
All of \citep{aravkin-baraldi-orban-2022,aravkin-baraldi-orban-2024,leconte-orban-2025} assume uniformly bounded second-order information in the model of \(f\).

\citet{leconte-orban-2023-2} revisit the trust-region method of \citep{aravkin-baraldi-orban-2022} but allow for unbounded model Hessians.
They establish global convergence and a worst-case complexity bound of \(O(\epsilon^{-2/(1-p)})\) provided \(\|B_k\| = O(|\mathcal{S}_k|^p)\) with \(0 \leq p < 1\).
To the best of our knowledge, they were the first to use that assumption and to obtain a complexity bound in the presence of unbounded model Hessians.
Unfortunately, their analysis does not generalize to \(p = 1\).
Under additional assumptions detailed in Section~5, they also prove the existence of a subsequence of iterates that converges to a first-order stationary point of~\eqref{eq:nlp}.
This result allows us to derive a similar convergence guarantee for both R2N and R2DH; nevertheless, their analysis still requires Lipschitz continuity of \(\nabla f\) in a neighborhood of the iterates.

Potentially unbounded model Hessians are a relevant assumption in several contexts, including quasi-Newton methods.
\citet[\S\(8.4.1.2\)]{conn-gould-toint-2000} show that the SR1 approximation satisfies \(\|B_k\| = O(|\mathcal{S}_k|)\), and a similar bound for BFGS when $f$ is convex.
\citet{powell-2010} establishes a similar bound for his PSB update.
Even though it is not currently known whether those bounds are tight, the case \(p = 1\) covers them.

\citet{diouane-habiboullah-orban-2024} generalized the results of \citep{leconte-orban-2023-2} to \(p = 1\) and provided tighter complexity constants when \(0 \leq p < 1\) in the context of trust-region methods for smooth optimization, i.e., \(h = 0\).
Our complexity analysis draws from \citep{diouane-habiboullah-orban-2024,leconte-orban-2023-2}.

\subsection*{Notation}

Unless otherwise noted, if \(x\) is a vector, \(\|x\|\) denotes its Euclidean norm and if \(A\) is a matrix, \(\|A\|\) denotes its spectral norm.
For positive sequences $\{a_k\}$ and $\{b_k\}$, we say that $a_k = o(b_k)$ if and only if $\limsup_k a_k / b_k = 0$.
The cardinality of a finite set $\mathcal{A}$ is denoted $|\mathcal{A}|$.
We denote \(\N_0\) the set of positive integers.

\section{Background}%
\label{sec:background}

We recall relevant concepts of variational analysis, e.g., \citep{rockafellar-wets-1998}.

The domain of \(h\) is \(\dom h := \{ x \in \R^n \mid h(x) < \infty\}\).
Because \(h\) is proper, \(\dom h \neq \varnothing\).
If \(P : \R^n \to \R^n\) is a set-valued function, \(\dom P = \{ x \in \R^n \mid P(x) \neq \varnothing \}\).

\begin{definition}%
  (Limiting subdifferential)
  Consider \(\phi : \R^n \to \widebar{\R}\) and \(\bar{x} \in \R^n\) such that \(\phi(\bar{x}) < +\infty\).
  We say that \(v \in \R^n\) is a regular subgradient of \(\phi\) at \(\bar{x}\) if
  \[
    \liminf_{x \to \bar{x}} \, \frac{\phi(x) - \phi(\bar{x}) - v^T (x - \bar{x})}{\|x - \bar{x}\|} \geq 0.
  \]
  The set \(\hat{\partial} \phi(\bar{x})\) of all regular subgradients of \(\phi\) at \(\bar{x}\) is called the Fréchet subdifferential.

  The limiting subdifferential of \(\phi\) at \(\bar{x}\) is the set \(\partial \phi(\bar{x})\) of all \(v \in \R^n\) such that there is \(\{x_k\} \to \bar{x}\) with \(\{\phi(x_k)\} \to \phi(\bar{x})\) and \(\{v_k\} \to v\) with \(v_k \in \hat{\partial} \phi(x_k)\) for all \(k\).
\end{definition}

If \(\phi = f + h\) with \(f\) continuously differentiable and \(h\) lower semi-continuous, then \(\partial \phi(x) = \nabla f(x) + \partial h(x)\) \citep[Theorem~\(10.1\)]{rockafellar-wets-1998}.

\begin{definition}%
  \label{def:prox}
  (Proximal Operator)
  Let \(h : \R^n \to \R \cup \{+\infty\}\) be proper and lower semi-continuous.
  The proximal operator of \(h\) with step length \(\nu > 0\) is
  \[
    \prox{\nu h}(x) := \argmin{y} \ h(y) + \tfrac{1}{2} \nu^{-1}\|y - x\|^2.
  \]
\end{definition}

Without further assumptions on \(h\), the proximal operator might be empty, or contain one or more elements.

By \citep[Exercise~\(8.8c\)]{rockafellar-wets-1998}, \(\bar{x}\) is first-order stationary for~\eqref{eq:nlp} if \(0 \in \nabla f(\bar{x}) + \partial h(\bar{x})\).

\begin{definition}%
  \label{def:lim-sup-osc}
  (Outer limit)
  Let \(\mathcal{K} \subseteq \N\) and a sequence \(\{S_k\}_{k \in \mathcal{K}}\) of subsets of \(\R^n\).
  The outer limit of \(\{S_k\}_{k \in \mathcal{K}}\) is
  \[
    \limsup_{k \in \mathcal{K}} S_k := \{ v \in \R^n \mid \exists \, \mathcal{K'} \subseteq \mathcal{K}, \, \exists \, v_k \in S_k \text{ for all } k \in \mathcal{K'} \text{ with } \lim_{k \in \mathcal{K'}} v_k = v \}.
  \]
\end{definition}

Let \(\mathcal{K} \subseteq \N\) and \(\{x_k\}_{k \in \mathcal{K}}\) such that \(\lim_{k \in \mathcal{K}} x_k = \bar{x}\) and \(\lim_{k \in \mathcal{K}} h(x_k) = h(\bar{x})\).
By \citep[Proposition~8.7]{rockafellar-wets-1998}, \(\limsup_{k \in \mathcal{K}} \partial h(x_k) \subseteq \partial h(\bar{x})\).

\section{Models}%
\label{sec:models}

For \(\sigma \geq 0\), \(x \in \R^n\), and \(B(x) = {B(x)}^T \in \R^{n \times n}\), consider the models
\begin{subequations}
  \begin{align}%
    \label{eq:def-phi}
    \varphi(s; x)   & := f(x) + \nabla {f(x)}^T s + \tfrac{1}{2} s^T B(x) s        \\
    \label{eq:def-psi}
    \psi(s; x)      & \phantom{:}\approx h(x + s)                                  \\
    \label{eq:def-m}
    m(s; x, \sigma) & := \varphi(s; x) + \tfrac{1}{2} \sigma \|s\|^2 + \psi(s; x).
  \end{align}
\end{subequations}
Note that~\eqref{eq:def-m} represents a regularized second-order model of the objective of~\eqref{eq:nlp}, where \(f\) and \(h\) are modeled separately.
More details, on the use of such model to solve~\eqref{eq:nlp}, will be given in \Cref{sec:algorithm}.
By construction, $\varphi(0;x)= f(x)$  and $\nabla \varphi(0;x)= \nabla f(x)$.
We make the following assumption on~\eqref{eq:def-psi}.

\begin{modelassumption}%
  \label{asm:psi}
  For any \(x \in \R^n\), \(\psi(\cdot; x)\) is proper, lower semi continuous and prox-bounded with threshold \(\lambda_x \in \R_+ \cup \{+\infty\}\) \citep[Definition~\(1.23\)]{rockafellar-wets-1998}.
  In addition, \(\psi(0; x) = h(x)\), and \(\partial \psi(0; x) = \partial h(x)\).
\end{modelassumption}
We make the following additional assumption and say that \(\{\psi(\cdot; x)\}\) is \emph{uniformly prox-bounded}.

\begin{modelassumption}%
  \label{asm:unif-prox-bounded}
  There is \(\lambda \in \R_+ \cup \{+\infty\}\) such that \(\lambda_x \geq \lambda\) for all \(x \in \R^n\).
\end{modelassumption}
\Cref{asm:unif-prox-bounded} is satisfied if \(h\) itself is prox-bounded and we select \(\psi(s; x) := h(x + s)\) for all \(x\).
Let
\begin{subequations}
  \begin{align}%
    p(x, \sigma) & := \min_s \ m(s; x, \sigma) \leq m(0; x, \sigma) = f(x) + h(x)
    \label{eq:def-p}%
    \\
    P(x, \sigma) & := \argmin{s} \ m(s; x, \sigma),
    \label{eq:def-P}
  \end{align}
\end{subequations}
be the value function and the set of minimizers of~\eqref{eq:def-m}, respectively.

For \(x \in \R^n\), $s \in P(x, \sigma) \Longrightarrow 0 \in \nabla \varphi(s; x) + \sigma s + \partial \psi(s; x)$.
Our first result states properties of the domain of \(p\) and \(P\) as given in~\eqref{eq:def-p} and~\eqref{eq:def-P}.

\begin{lemma}%
  \label{lem:domains}
  Let \Cref{asm:psi} be satisfied and \(B(x) = B(x)^T\) for all $x \in \R^n$.
  Then, \(\dom p = \R^n \times \R\).
  In addition, if \Cref{asm:unif-prox-bounded} holds, \(\dom P \supseteq \{ (x, \sigma) \mid \sigma > \max(\lambda^{-1} -\lambda_{\min}(B(x)), \, \lambda^{-1}) \}\), where \(\lambda_{\min}(B(x))\) is the smallest eigenvalue of \(B(x)\).
\end{lemma}

\begin{proof}
  By definition of the domain and \Cref{asm:psi},
  \begin{align*}
    \dom p  = \{ (x, \sigma)  \mid \inf_s m(s; x, \sigma) < +\infty \} & = \{ (x, \sigma)  \mid \exists s \ m(s; x, \sigma) < +\infty \}
    \\
                                                                       & = \{ (x, \sigma)  \mid \exists s \ \psi(s; x) < +\infty \} = \R^n \times \R,
  \end{align*}
  because \(\psi(\cdot; x)\) is proper.
  Moreover,
  \begin{equation*}
    \dom P = \{ (x, \sigma) \mid \exists s(x, \sigma) \in \R^n, \ m(s(x, \sigma); x, \sigma) = \inf_s m(s; x, \sigma) \}.
  \end{equation*}
  Write
  \[
    m(s; x, \sigma) = \varphi(s; x) + \tfrac{1}{2} (\sigma - \lambda^{-1}) \|s\|^2 + \psi(s; x) + \tfrac{1}{2} \lambda^{-1} \|s\|^2.
  \]
  By \Cref{asm:unif-prox-bounded} and \citep[Exercise~\(1.24(c)\)]{rockafellar-wets-1998}, there is \(b \in \R\) such that \(\psi(s; x) + \tfrac{1}{2} \lambda^{-1} \|s\|^2 \geq b \) for all \(s \in \R^n\).
  Let \(a \in \R\).
  The above and~\eqref{eq:def-phi} imply that the level set \(\{s \in \R^n \mid m(s; x, \sigma) \leq a \}\) is contained in
  \[
    \{ s\in \R^n \mid \nabla {f(x)}^T s + \tfrac{1}{2} s^T (B(x) + (\sigma - \lambda^{-1}) I) s \leq a - b - f(x) \},
  \]
  which is a bounded set for \(\sigma > \lambda^{-1} - \lambda_{\min}(B(x))\), i.e., \(m(\cdot; x, \sigma)\) is level-bounded.
  Thus, \citep[Theorem~\(1.9\)]{rockafellar-wets-1998} implies that \(\inf_s m(s; x, \sigma)\) is attained, i.e., that \(P(x, \sigma) \neq \varnothing\).
\end{proof}

In \Cref{lem:domains}, \(\dom P = \{ (x, \sigma) \mid \sigma > \max(\lambda^{-1} -\lambda_{\min}(B(x)), \, \lambda^{-1}) \}\) does not hold in general.
Consider for example a situation where \(\psi(s; x)\) is bounded below for all \(x \in \R^n\), i.e., each \(\lambda_x = +\infty\).
We can choose \(\lambda = +\infty\).
Assume also that, for a given \(x \in \R^n\), \(\varphi(s; x) = 0\) for all \(s\), and \(\psi(s; x)\) level-bounded.
Then, \(\lambda_{\min}(B(x)) = 0\), and for \(\sigma = 0 = \lambda^{-1}\), \(m(s; x, \sigma) = \psi(s; x)\).
Therefore, \(P(x, \sigma) \neq \varnothing\).

For a given \(s \in P(x, \sigma)\), we define
\begin{equation}%
  \label{eq:def-xi}
  \xi(s; x, \sigma) := f(x) + h(x) - (\varphi(s; x) + \psi(s; x)).
\end{equation}

The next result relates~\eqref{eq:def-xi} to first-order stationary for~\eqref{eq:nlp} and~\eqref{eq:def-m}.

\begin{lemma}%
  \label{lem:stationarity}
  Let \Cref{asm:psi} be satisfied, and \(x \in \R^n\) and \(\sigma \geq 0\) be given.
  Then, for \(s \in P(x, \sigma)\), \(\xi(s; x, \sigma) = 0 \ \Longrightarrow \ s = 0 \ \Longrightarrow \ x\) is first-order stationary for~\eqref{eq:nlp}.
\end{lemma}

\begin{proof}
  For \(s \in P(x, \sigma)\), if \(\xi(s; x, \sigma) = 0\), then
  \begin{equation*}
    0 = \xi(s; x, \sigma) = f(x) + h(x) - ( \varphi(s; x) + \psi(s; x) ) \ge \frac{1}{2} \sigma \|s\|^2,
  \end{equation*}
  which implies that \(s = 0\), and therefore \(0 \in P(x, \sigma)\).
  Therefore, \(0 \in \partial m(0; x, \sigma) = \nabla \varphi(0; x) + \partial \psi(0; x) = \nabla f(x) + \partial h(x)\), and \(x\) is first-order stationary for~\eqref{eq:nlp}.
\end{proof}

The following proposition states some properties of~\eqref{eq:def-p} and~\eqref{eq:def-P}.

\begin{proposition}%
  \label{lem:prop-p}
  Let \Cref{asm:psi,asm:unif-prox-bounded} be satisfied.
  Assume also that \(\nabla f\) is bounded over \(\R^n\).
  Let \(\epsilon > 0\).
  Then,
  \begin{enumerate}
    \item at any \((x, \sigma)\) such that \(\sigma \geq \lambda^{-1} - \lambda_{\min}(B(x)) + \epsilon\), \(p\) is finite and lsc, and \(P(x, \sigma)\) is nonempty and compact;
    \item%
      \label{itm:p-cont}
      if \(\{(x_k, \sigma_k)\} \to (\bar{x}, \bar{\sigma})\) with \(\sigma_k \geq \lambda^{-1} - \lambda_{\min}(B(x_k)) + \epsilon\) for all \(k\) in such a way that \(\{p(x_k, \sigma_k)\} \to p(\bar{x}, \bar{\sigma})\), and for each \(k\), \(s_k \in P(x_k, \sigma_k)\), then \(\{s_k\}\) is bounded and all its limit points are in \(P(\bar{x}, \bar{\sigma})\);
    \item for any \(x \in \R^n\), \(p(\bar{x}, \cdot)\) is continuous at any \(\bar{\sigma} \geq \lambda^{-1} - \lambda_{\min}(B(\bar{x})) + \epsilon\) and \(\{p(x_k, \sigma_k)\} \to p(\bar{x}, \bar{\sigma})\) holds in part~\ref{itm:p-cont}.
  \end{enumerate}
\end{proposition}

\begin{proof}
  The proof consists in establishing that~\eqref{eq:def-m} is level-bounded in \(s\) locally uniformly in \((x, \sigma)\) \citep[Definition~\(1.16\)]{rockafellar-wets-1998} for \(\sigma \geq \lambda^{-1} - \lambda_{\min}(B(x)) + \epsilon\) and applying \citep[Theorem~\(1.17\)]{rockafellar-wets-1998}.
  It is nearly identical to that of \citep[Proposition~\(3.2\)]{aravkin-baraldi-orban-2024} and is omitted.
\end{proof}

Even though model~\eqref{eq:def-m} is natural for incorporating second-order information, it is generally difficult to compute an exact minimizer of it.
We proceed as \citet{aravkin-baraldi-orban-2022,aravkin-baraldi-orban-2024} and consider a simpler first-order model that will allow us to define an implementable stationary measure, to set minimal requirements steps computed in the course of the iterations of the algorithm of \Cref{sec:algorithm}, and to derive convergence properties.
This first-order model generalizes the concept of \emph{Cauchy point} (``cp'') when solving~\eqref{eq:nlp}.
For fixed \(\nu > 0\) and \(x \in \R^n\), define
\begin{subequations}%
  \label{eq:cp-models}
  \begin{align}
    \varphicp(s; x)                 & := f(x) + \nabla {f(x)}^T s                                                                    \\
    m_{\textup{cp}}(s; x, \nu^{-1}) & := \varphicp(s; x) + \tfrac{1}{2} \nu^{-1} \|s\|^2 + \psi(s; x)
    \label{eq:def-m-cp}%
    \\
    p_{\textup{cp}}(x, \nu^{-1})    & := \min_s \ m_{\textup{cp}}(s; x, \nu^{-1}) \leq m_{\textup{cp}}(0; x, \nu^{-1}) = f(x) + h(x)
    \label{eq:def-p-cp}%
    \\
    P_{\textup{cp}}(x, \nu^{-1})    & := \argmin{s} \ m_{\textup{cp}}(s; x, \nu^{-1})
    \label{eq:def-P-cp}%
    \\
    \xicp(\scp; x, \nu^{-1})        & := f(x) + h(x) - (\varphicp(\scp; x) + \psi(\scp; x)),
    \label{eq:def-xi1}
  \end{align}
\end{subequations}
where \(\scp \in P_{\textup{cp}}(x, \nu^{-1})\).
By \citep[Lemma~\(2\)]{bolte-sabach-teboulle-2014},
\begin{equation}%
  \label{eq:xi-decrease}
  \xicp(\scp; x, \nu^{-1}) \geq \tfrac{1}{2} \nu^{-1} \|\scp(x, \nu^{-1})\|^2 \ge 0.
\end{equation}

In the smooth case, i.e., $h = 0$ and \(\psi = 0\), \(\scp = - \nu \nabla f(x)\), so that
\[
  \xicp(\scp; x, \nu^{-1}) = f(x) - (f(x) + \nabla f(x)^T \scp) = \nu \| \nabla f(x) \|^2,
\]
which suggests \(\nu^{-1/2} {\xicp(\scp; x, \nu^{-1}) }^{1/2}\) as a stationarity measure that generalizes the norm of the gradient to the nonsmooth setting.

Furthermore, this choice can be naturally interpreted in terms of the subdifferential of the model \(m_{\textup{cp}}\) \cref{eq:def-m-cp} at \(\scp\).
Specifically, from~\eqref{eq:xi-decrease}, we have
\begin{equation}
  \label{eq:xi-norm}
  \sqrt{\nu^{-1} \xi_{\textup{cp}}(s_{\textup{cp}}; x, \nu^{-1})} \geq \tfrac{1}{\sqrt{2}} \nu^{-1} \|s_{\textup{cp}}\|.
\end{equation}
Moreover, by definition of \(s_{\textup{cp}} \in P_{\textup{cp}}(x, \nu^{-1})\), the first-order optimality condition yields
\[
  -\nu^{-1} s_{\textup{cp}} \in \nabla f(x) + \partial \psi(s_{\textup{cp}}; x) = \partial m_{\textup{cp}}(s_{\textup{cp}}; x, \nu^{-1}).
\]
Hence,
\[
  \operatorname{dist}(0, \partial m_{\textup{cp}}(s_{\textup{cp}}; x, \nu^{-1})) \leq \nu^{-1} \|s_{\textup{cp}}\| \leq \sqrt{2 \nu^{-1} \xi_{\textup{cp}}(x, \nu^{-1})},
\]
which directly links our stationarity measure to the distance to the subdifferential of the first-order model \(m_{\textup{cp}}\) at \(\scp\).

It is worth noting that other authors, such as \citet{kanzow-mehlitz-2022}, adopt a different stationarity measure by directly considering \(\nu^{-1} \|s_{\textup{cp}}\|\).

The next results establish corresponding properties of \(p_{\textup{cp}}\) and \(P_{\textup{cp}}\).
The proofs are similar to those of \Cref{lem:domains,lem:stationarity,lem:prop-p} and are omitted.

\begin{lemma}%
  \label{lem:domains-cp}
  Let \Cref{asm:psi} be satisfied.
  Then, \(\dom p_{\textup{cp}} = \R^n \times \R\).
  If \Cref{asm:unif-prox-bounded} holds, \(\dom P_{\textup{cp}} \supseteq \{(x, \nu^{-1}) \mid \nu > \max(\lambda^{-1} - \lambda_{\min}(B(x)), \, \lambda^{-1}) \}\).
\end{lemma}

The next result characterizes first-order stationarity for~\eqref{eq:nlp}.

\begin{lemma}%
  \label{lem:stationarity-cp}
  Let \Cref{asm:psi} be satisfied and \(\nu > 0\).
  Then, for \(\scp \in P_{\textup{cp}}(x, \nu^{-1})\), \(\xicp(\scp; x, \nu^{-1}) = 0 \ \Longrightarrow \ \scp = 0 \ \Longrightarrow \ x\) is first-order stationary for~\eqref{eq:nlp}.
\end{lemma}

The following result states properties of~\eqref{eq:def-p-cp} and~\eqref{eq:def-P-cp}.

\begin{proposition}%
  \label{lem:prop-p-cp}
  Let \Cref{asm:psi,asm:unif-prox-bounded} be satisfied and \(\nabla f(x)\) be bounded over \(\R^n\).
  Let \(\epsilon > 0\).
  Then,
  \begin{enumerate}
    \item at any \((x, \nu^{-1})\) with \(\nu^{-1} \geq \lambda^{-1} + \epsilon\), \(p_{\textup{cp}}\) is finite and lsc, and \(P_{\textup{cp}}(x, \nu^{-1})\) is nonempty and compact;
    \item%
      \label{itm:p-cont-cp}
      if \(\{(x_k, \nu_k^{-1})\} \to (\bar{x}, \bar{\nu}^{-1})\) with \(\nu_k^{-1} \geq \lambda^{-1} + \epsilon\) for all \(k\) in such a way that \(\{p_{\textup{cp}}(x_k, \nu_k^{-1})\} \to p_{\textup{cp}}(\bar{x}, \bar{\nu}^{-1})\), and for each \(k\), \(s_k \in P_{\textup{cp}}(x_k, \nu_k^{-1})\), then \(\{s_k\}\) is bounded and all its limit points are in \(P_{\textup{cp}}(\bar{x}, \bar{\nu}^{-1})\);
    \item for any \(\bar{x} \in \R^n\) and any \(\bar{\nu}^{-1} \geq \lambda^{-1} + \epsilon\), \(p_{\textup{cp}}(\bar{x}, \cdot)\) is continuous at \(\bar{\nu}\) and \(\{p_{\textup{cp}}(x_k, \nu_k^{-1})\} \to p_{\textup{cp}}(\bar{x}, \bar{\nu}^{-1})\) holds in part~\ref{itm:p-cont-cp}.
  \end{enumerate}
\end{proposition}

The main idea of the algorithm proposed in \Cref{sec:algorithm} is that~\eqref{eq:def-m} is approximately minimized at each iteration.
In order to establish convergence, the step \(s\) thus computed is required to satisfy \emph{Cauchy decrease}, which we define as in \citep{aravkin-baraldi-orban-2024,aravkin-baraldi-orban-2022}:
\begin{equation}%
  \label{eq:cauchy-decrease}
  \varphi(0; x) + \psi(0; x) - (\varphi(s; x) + \psi(s; x)) \geq (1-\theta_1) \xicp(\scp; x, \nu^{-1}) ,
\end{equation}
for a preset value of \(\theta_1 \in (0, \, 1)\).
In other words, \(s\) must result in a decrease in \(\varphi(\cdot; x) + \psi(\cdot; x)\) that is at least a fraction of the decrease of the Cauchy model \(\varphicp(\cdot; x) + \psi(\cdot; x)\) obtained with the Cauchy step \(\scp\) and a well-chosen step length \(\nu\).

The following result parallels \citep[Proposition~4]{leconte-orban-2023-2} and establishes that if a step \(s\) reduces~\eqref{eq:def-m} at least as much as \(\scp\) does, Cauchy decrease holds.
This observation is important because the first step of the proximal-gradient method from \(s = 0\) applied to~\eqref{eq:def-m-cp} and to~\eqref{eq:def-m} with step length \(\nu\) is the same, and that step is precisely \(\scp\).
Therefore, a step \(s\) may be obtained by continuing the proximal-gradient iterations on~\eqref{eq:def-m} from \(\scp\).

\begin{proposition}%
  \label{prop:step_asm}
  Let \Cref{asm:psi} be satisfied.
  Let $x \in \R^n$, $\theta_1 \in (0, \, 1)$, $\sigma > 0$ and let \(\scp\) be computed with $\nu = \theta_1 / (\|B(x)\| + \sigma)$.
  Assume $s\in \R^n$ is such that $m(s; x, \sigma) \le m(\scp; x, \sigma)$.
  Then, \(s\) satisfies~\eqref{eq:cauchy-decrease}.
\end{proposition}

\begin{proof}
  Let $x \in \R^n$, $\sigma>0$, and $s\in \R^n$, such that $m(s; x, \sigma) \le m(\scp; x,\sigma)$.
  Then,
  \begin{align*}
    \varphi(s; x) + \psi(s; x) + \tfrac{1}{2} \sigma \|s\|^2 & \leq \varphi(\scp; x) + \psi(\scp; x) + \tfrac{1}{2} \sigma \|\scp\|^2
    \\
                                                             & = \varphicp(\scp; x) + \psi(\scp; x) + \tfrac{1}{2} \scp^T B(x) \scp + \tfrac{1}{2} \sigma \|\scp\|^2.
  \end{align*}
  The Cauchy-Schwarz inequality \(\scp^T B(x) \scp \leq \|B(x)\| \|\scp\|^2\), the identity \(\varphi(0; x) = \varphicp(0; x)\)  and~\eqref{eq:xi-decrease} yield
  \begin{align*}
    (\varphi + \psi)(0; x) - (\varphi + \psi)(s; x) & \geq \xicp(\scp; x, \nu^{-1}) - \tfrac{1}{2} (\|B(x)\| + \sigma) \|\scp\|^2 + \tfrac{1}{2} \sigma \|s\|^2
    \\
                                                    & \geq \xicp(\scp; x, \nu^{-1}) - \tfrac{1}{2} (\|B(x)\| + \sigma) \|\scp\|^2
    \\
                                                    & \geq \xicp(\scp; x, \nu^{-1}) - (\|B(x)\| + \sigma) \nu \xicp(\scp; x, \nu^{-1})
    \\
                                                    & = (1-\theta_1) \xicp(\scp; x, \nu^{-1}) .
  \end{align*}
\end{proof}

Computing \(\|B(x)\|\) in the spectral norm comes at a cost.
However, as we now illustrate, an inexact computation is sufficient in order to ensure~\eqref{eq:cauchy-decrease}.
Assume that we are able to compute \(\beta(x) \approx \|B(x)\|\) such that \(\beta(x) \geq \mu \|B(x)\|\) for \(0 < \mu < 1\), and set \(\nu = \theta_1 / (\beta(x) + \sigma)\).
The proof of \Cref{prop:step_asm} continues to apply unchanged until the very last line, which becomes
\begin{align*}
  \varphi(0; x) + \psi(0; x) - (\varphi(s; x) + \psi(s; x)) & \geq (1 -  (\|B(x)\| + \sigma) \nu) \xicp(\scp; x, \nu^{-1})
  \\
                                                            & = \left( 1 - \theta_1 \frac{\|B(x)\| + \sigma}{\beta(x) + \sigma} \right) \xicp(\scp; x, \nu^{-1}) .
\end{align*}
If \(\beta(x) \leq \|B(x)\|\), \((\|B(x)\| + \sigma) / (\beta(x) + \sigma) \leq \|B(x)\| / \beta(x) \leq 1 / \mu\), so that
\[
  \left( 1 - \theta_1 \frac{\|B(x)\| + \sigma}{\beta(x) + \sigma} \right) \xicp(\scp; x, \nu^{-1}) \geq (1 - \theta_1/\mu) \xicp(\scp; x, \nu^{-1}) .
\]
Thus, as long as \(\theta_1 < \mu\),~\eqref{eq:cauchy-decrease} is satisfied with \(\theta_1\) replaced with \(\theta_1 / \mu\).

If, on the other hand, \(\beta(x) \geq \|B(x)\|\), then \((\|B(x)\| + \sigma) / (\beta(x) + \sigma) \leq 1\), and
\[
  \left( 1 - \theta_1 \frac{\|B(x)\| + \sigma}{\beta(x) + \sigma} \right) \xicp(\scp; x, \nu^{-1}) \geq (1 - \theta_1) \xicp(\scp; x, \nu^{-1}) ,
\]
and~\eqref{eq:cauchy-decrease} holds unchanged.

The above observation also allows us to replace \(\|B(x)\|\) in the denominator of \(\nu\) with, e.g., \(\|B(x)\|_1\), \(\|B(x)\|_{\infty}\) or \(\|B(x)\|_F\) if \(B(x)\) is available as an explicit matrix, or indeed with any other norm of \(B(x)\).

\section{A modified quasi-Newton method for nonsmooth optimization}%
\label{sec:algorithm}

We are in position to describe a modified quasi-Newton method to solve~\eqref{eq:nlp} named R2N.
By contrast with trust-region-based approaches \citep{aravkin-baraldi-orban-2022,leconte-orban-2025}, proximal operators are easier to evaluate in the R2N subproblem as they do not include a trust-region indicator.

At iteration \(k\), we choose a step length \(\nu_k > 0\) based on the regularization parameter \(\sigma_k > 0\) and the norm of the model Hessian \(B(x_k)\) at the current iterate \(x_k \in \R^n\) as in \Cref{prop:step_asm}.
We then compute the Cauchy step \(\skcp\) as a minimizer of~\eqref{eq:def-m-cp}.
A step \(s_k\) is subsequently computed that satisfies the assumptions of \Cref{prop:step_asm}.

The rest of the algorithm is standard.
The decrease in \(f + h\) at $x_k + s_k$ is compared to the decrease predicted by the model.
If both are in sufficient agreement, $x_k + s_k$ becomes the new iterate, and \(\sigma_k\) is possibly reduced.
If the model turns out to predict poorly the actual decrease, the trial point is rejected and \(\sigma_k\) is increased.
\Cref{alg:R2N} states the whole procedure.

The interaction between \(\sigma_k\) and the unknown threshold \(\lambda_{x_k}\) works as in \citep[Algorithm~\(6.1\)]{aravkin-baraldi-orban-2022} and \citep{aravkin-baraldi-orban-2024}.
If \(\sigma_k \leq \lambda_{x_k}^{-1}\), \(\psi(s_k; x_k) = -\infty\), and according to the rules of extended arithmetic, which state that \(\pm\infty \cdot 0 = 0 \cdot (\pm\infty) = (\pm\infty)/(\pm\infty) := 0\) \citep{rockafellar-wets-1998}, \(\rho_k = 0\).
Consequently, \(s_k\) will be rejected at Line~\ref{alg:R2N:step-accept}, and \(\sigma_{k+1}\) will be set larger than \(\sigma_k\) at Line~\ref{alg:R2N:step-update}.
After a finite number of such increments, \(\sigma_k\) will surpass \(\lambda_{x_k}^{-1}\), resulting in a step with finite \(\psi(s_k; x_k)\).
In effect, \Cref{asm:unif-prox-bounded} is only required to hold at the iterates generated by the algorithm.

\begin{algorithm}[ht]%
  \caption[caption]{%
    \label{alg:R2N}
    R2N\@: A proximal modified Quasi-Newton method.
  }
  \begin{algorithmic}[1]%
    \State Choose constants  \(0 < \theta_1 < 1 < \theta_2 \), \(0 < \eta_1 \leq \eta_2 < 1\) and \(0 < \gamma_3 \leq 1 < \gamma_1 \leq \gamma_2\).
    \State Choose \(\sigma_0 > 0\) and \(x_0 \in \R^n\) where \(h\) is finite.
    \For{\(k = 0, 1, \dots\)}
    \State%
    \label{alg:R2N:Bk}
    Choose \(B_k:=B(x_k) \in \R^{n \times n}\) such that \(B_k = B_k^T \).
    \State%
    \label{alg:R2N:step-nuk}
    Compute \(\nu_k:= \theta_1 / (\|B_k\| + \sigma_k)\).
    \State%
    \label{alg:R2N:step-switch}
    Compute \(\skcp \in \argmin{s} \ m_{\textup{cp}}(s; x_k, \nu_k^{-1})\) and \(\xicp(\skcp, x_k, \nu_k^{-1})\) as defined in~\eqref{eq:def-xi1}.
    \State%
    \label{alg:R2N:step-computation}
    Compute a step \(s_k\) such that \(m(s_k; x_k, \sigma_k)\le m(\skcp; x_k, \sigma_k) \).
    \State%
    \label{alg:R2N:step-comparison}
    If \(\|s_k\| > \theta_2 \; \|\skcp\| \), reset \(s_k = \skcp\).
    \State%
    \label{alg:R2N:step-rhok}
    Compute the ratio
    \[
      \rho_k :=
      \frac{
        f(x_k) + h(x_k) - (f(x_k + s_k) + h(x_k + s_k))
      }{
        \varphi(0; x_k) + \psi(0; x_k) - (\varphi(s_k; x_k) + \psi(s_k; x_k))
      }.
    \]
    \State%
    \label{alg:R2N:step-accept}%
    If \(\rho_k \geq \eta_1\), set \(x_{k+1} = x_k + s_k\).
    Otherwise, set \(x_{k+1} = x_k\).
    \State%
    \label{alg:R2N:step-update}%
    Update the regularization parameter according to
    \[
      \sigma_{k+1} \in
      \begin{cases}
        \begin{aligned}
           & [\gamma_3 \sigma_k, \, \sigma_k]          &  & \text{ if } \rho_k \geq \eta_2,          &  & \quad \text{very successful iteration} \\
           & [\sigma_k, \, \gamma_1 \sigma_k]          &  & \text{ if } \eta_1 \leq \rho_k < \eta_2, &  & \quad \text{successful iteration}      \\
           & [\gamma_1 \sigma_k, \, \gamma_2 \sigma_k] &  & \text{ if } \rho_k < \eta_1.             &  & \quad \text{unsuccessful iteration}
        \end{aligned}
      \end{cases}
    \]
    \EndFor
  \end{algorithmic}
\end{algorithm}

Importantly, R2N does not require \(B_k \succeq 0\), which may be useful in practice in order to capture natural problem curvature.
In addition, we allow \(\{B_k\}\) to be unbounded.
In \Cref{sec:convergence}, we establish convergence provided it does not diverge too fast, using an assumption similar to that used in trust-region methods \citep[\S\(8.4\)]{conn-gould-toint-2000}.
In \Cref{sec:evaluation-complexity}, we study the effect of using  different bounds on \(\{\|B_k\|\}\) on worst-case evaluation complexity.
The complexity results are obtained by adapting results from \citep{diouane-habiboullah-orban-2024} (in the context of trust-region methods for smooth optimization) to R2N.

Our main working assumption is the following.

\begin{problemassumption}%
  \label{asm:g-continuous}
  The function \(f\) is continuously differentiable over the set \( \{x \in \R^n \mid (f+h)(x) \leq (f+h)(x_0)\} \) and \(h\) is proper and lower semi-continuous.
\end{problemassumption}

\Cref{asm:g-continuous} is very mild as one does not require boundedness nor Lipschitz continuity of \(f\) or \(\nabla f\), in contrast with \citep[Problem Assumption~\(4.1\)]{aravkin-baraldi-orban-2024} or the assumptions of \citet{kanzow-lechner-2024}.
For instance, our analysis includes cases where \(f\) is continuously differentiable, but whose gradient is not locally Lipschitz continuous at \(x = 0\), e.g., \(f(x) = |x|^\frac{3}{2}\).

In the next sections, we derive convergence and worst-case complexity analysis for \Cref{alg:R2N}.
We will repeatedly use the notation
\begin{subequations}%
  \label{eq:S-U-sets}
  \begin{alignat}{2}%
    \mathcal{S}   & := \{ i \in \N \mid \rho_i \geq \eta_1 \}              &  & \qquad \text{(all successful iterations)}
    \\
    \mathcal{S}_k & := \{ i \in \mathcal{S} \mid i \le k \}                &  & \qquad \text{(successful iterations until iteration \(k\))}
    \\
    \mathcal{U}   & := \{ i \in \N \mid \rho_i < \eta_1 \}                 &  & \qquad \text{(all unsuccessful iterations)}
    \\
    \mathcal{U}_k & := \{i \in \N \mid i\not \in \mathcal{S}, \ i \le k \} &  & \qquad \text{(unsuccessful iterations until iteration \(k\))}.
  \end{alignat}
\end{subequations}

\section{Convergence analysis of \texorpdfstring{\Cref{alg:R2N}}{Algorithm 4.1}}%
\label{sec:convergence}

In this section, we investigate the convergence properties of \Cref{alg:R2N} under \Cref{asm:g-continuous}.
For notational convenience, we denote \(\xicp(\skcp, x_k, \nu_k^{-1})\) by \(\xikcp\).
We then show that \(\liminf_{k \to \infty} \nu_k^{-\frac{1}{2}} {\xikcp}^{\frac{1}{2}} = 0\).
We stress that the obtained convergence properties of \Cref{alg:R2N} are more general than those of \citep{aravkin-baraldi-orban-2024,kanzow-lechner-2024,kanzow-mehlitz-2022}, and do not require boundedness of the model Hessians nor (local) Lipschitz continuity of \(\nabla f\).

We first establish lower bounds on $\xikcp$ in terms of $\|s_k\|$.
\begin{lemma}%
  \label{lem:lm-bound-Seps}
  For all \(k \in \N\),
  \begin{equation}
    \label{eq:lm-bound-Seps2}
    \xikcp \geq \tfrac{1}{2 \theta_2^2} \nu_k^{-1} \|s_k\|^2.
  \end{equation}
  Additionally, for any $\alpha>0$,
  \begin{equation}
    \label{lem:xi-bound}
    \nu_k^{-\frac{1}{2}}{\xikcp}^{\frac{1}{2}} \geq \alpha  \quad \Rightarrow \quad   \xicp(x_k; \nu_k^{-1}) \geq \frac{\alpha}{\theta_2 \sqrt{2}} \|s_k\|.
  \end{equation}
\end{lemma}

\begin{proof}
  From \Cref{alg:R2N}, we have $\|s_k\| \le \theta_2 \|\skcp\|$.
  Hence,
  \begin{equation*}
    \label{ineq:xi_cp}
    \xikcp \ge \tfrac{1}{2} \nu_k^{-1} \|\skcp\|^2
    \ge \tfrac{1}{2\theta_2^2} \nu_k^{-1} \|s_k\|^2.
  \end{equation*}
  If $\nu_k^{-\frac{1}{2}}{\xikcp}^{\frac{1}{2}} \geq \alpha$,
  \begin{align*}
    \xikcp \geq \alpha \nu_k^\frac{1}{2} \xikcp^\frac{1}{2} \geq \alpha \nu_k^\frac{1}{2} {( \tfrac{1}{2 \theta_2^2} \nu_k^{-1} \|s_k\|^2 )}^\frac{1}{2} = \frac{\alpha}{\theta_2 \sqrt{2}} \|s_k\|.
  \end{align*}
\end{proof}

The next lemma shows that the convergence of \({\{x_k\}}_{k \in \N}\) holds if the objective is bounded below, the algorithm generates infinitely many successful iterations and the stationarity measure \(\nu_k^{-1/2} {\xikcp}^{1/2}\) is bounded away from zero.

\begin{lemma}%
  \label{lem:xk-convergennt}
  Assume that \Cref{alg:R2N} generates infinitely many successful iterations and that there is \({(f+h)}_{\text{low}} \in \mathbb{R}\) such that \((f+h)\left(x_k\right) \ge {(f+h)}_{\text{low}}\) for all \(k \in \mathbb{N}\).
  Additionally, assume, that there is \(\alpha > 0\) such that for all \(k \in \N\), \( \nu_k ^{-1/2}\xikcp^{1/2} \ge \alpha \).
  Then, \({\{x_k\}}_{k \in \N}\) is a Cauchy sequence, and hence converges.
\end{lemma}

\begin{proof}
  For all \(k \in \mathcal{S}\), using~\eqref{lem:xi-bound} from \Cref{lem:lm-bound-Seps}, we have
  \begin{align*}
    f(x_k) + h(x_k) - f(x_{k + 1}) - h(x_{k + 1}) & \geq \eta_1 (1-\theta_1) \xikcp
    \\
                                                  & \geq \frac{\eta_1 (1-\theta_1) \alpha}{\theta_2 \sqrt{2}} \|s_k\| = \frac{\eta_1 (1-\theta_1) \alpha}{\theta_2 \sqrt{2}} \|x_{k+1} - x_k\|.
  \end{align*}
  Summing over all successful iterations from 1 to \(k\), we obtain
  \begin{align*}
    f(x_0) + h(x_0) - {(f+h)}_{\text{low}} & \geq \sum_{j \in \mathcal{S}_k} f(x_j) + h(x_j) - f(x_{j+1}) - h(x_{j+1})
    \\
                                           & \geq \frac{\eta_1 (1 - \theta_1) \alpha }{\theta_2 \sqrt{2}} \sum_{j \in \mathcal{S}_k} \|x_{j+1} - x_j\|
    \\
                                           & = \frac{\eta_1 (1 - \theta_1) \alpha}{\theta_2 \sqrt{2}} \sum_{j=0}^k \left\|x_{j+1} - x_j\right\|.
  \end{align*}
  Thus, \( \sum_{j \in \N} \|x_{j+1} - x_j\| < +\infty\).
  Hence, \({\{x_k\}}_{k \in \N}\) is a Cauchy sequence, and converges.
\end{proof}

The following lemma shows that when \(\{x_k\}\) converges and \(\{\sigma_k\}\) diverges along common subsequences, the corresponding subsequence of \(\{s_k\}\) converges to zero.

\begin{lemma}%
  \label{lem:unbounded-sigma}
  Let \Cref{asm:g-continuous,asm:psi,asm:unif-prox-bounded} be satisfied and assume that there is an index set \(\mathcal{K} \subseteq \N\) such that
  \(\lim_{k \in \mathcal{K}} \sigma_k = + \infty \) and \(\{x_k\}_{k \in \mathcal{K}}\) is bounded.
  Then,
  \(\lim_{k \in \mathcal{K}} s_k = \lim_{k \in \mathcal{K}} \skcp = 0\).
\end{lemma}

\begin{proof}
  By contradiction, assume that there is an index set \(\mathcal{K}' \subseteq \mathcal{K}\) and \(\alpha>0\) such that \(\|\skcp\| \geq \alpha \) for all \(k \in \mathcal{K}'\).
  By definition of \(\skcp\) and \Cref{asm:psi}, \(f(x_k) + h(x_k)=m_{\textup{cp}}(0; x_k, \nu_k^{-1}) \ge m_{\textup{cp}}(\skcp; x_k, \nu_k^{-1}) \).
  Hence,
  \begin{align*}
    f(x_k) + h(x_k) & \geq \varphicp(\skcp; x_k) + \psi(\skcp; x_k) + \tfrac{1}{2} \nu_k^{-1} \|\skcp\|^2                               \\
                    & = f(x_k) + \nabla {f(x_k)}^T \skcp + \tfrac{1}{2\theta_1} (\|B_k\| + \sigma_k) \|\skcp\|^2   + \psi(\skcp; x_k)   \\
                    & \geq f(x_k) + \nabla {f(x_k)}^T \skcp + \tfrac{1}{2\theta_1} \sigma_k \|\skcp\|^2   + \psi(\skcp; x_k)            \\
                    & \geq f(x_k) - \|\nabla {f(x_k)}\| \|\skcp\| + \tfrac{1}{2} (\frac{\sigma_k}{\theta_1} - \lambda^{-1}) \|\skcp\|^2 \\
                    & \phantom{\geq} \ + \psi(\skcp; x_k) + \tfrac{1}{2} \lambda^{-1} \|\skcp\|^2.
  \end{align*}
  By \Cref{asm:unif-prox-bounded} and \citep[Exercise~\(1.24(c)\)]{rockafellar-wets-1998}, there is \(b_h \in \R\) such that \(\psi(s; x) + \tfrac{1}{2} \lambda^{-1} \|s\|^2 \ge b_h\) for all $s$ and $x$.
  Hence, for all sufficiently large \(k \in \mathcal{K}'\), \(\sigma_k > \lambda^{-1}\) and
  \begin{align}
    f(x_k) + h(x_k) & \geq f(x_k) - \|\nabla {f(x_k)}\| \|\skcp\| + \tfrac{1}{2} (\frac{\sigma_k}{\theta_1}  - \lambda^{-1}) \|\skcp\|^2 + b_h  \nonumber                      \\
                    & \geq f(x_k) - \|\nabla {f(x_k)}\| \|\skcp\| + \tfrac{1}{2} \alpha (\frac{\sigma_k}{\theta_1}  - \lambda^{-1}) \|\skcp\| + b_h                  \nonumber \\
                    & = f(x_k) + \left( \tfrac{1}{2} \alpha (\frac{\sigma_k}{\theta_1}  - \lambda^{-1}) - \|\nabla {f(x_k)}\| \right) \|\skcp\| + b_h.
    \label{eq:lm:bound:step}
  \end{align}
  Since \({\{x_k\}}_{k \in \mathcal{K}'}\) is bounded, so are \({\{f(x_k)\}}_{k \in \mathcal{K}'}\) and \({\{\nabla {f(x_k)}\}}_{k \in \mathcal{K}'}\) by \Cref{asm:g-continuous}.
  Let \(b_f := \min_{k \in \mathcal{K}'} f(x_k) > -\infty\) and \(b_{f'} = \max_{k \in \mathcal{K}'} \|\nabla {f(x_k)}\| < \infty\).
  Because \(\{f(x_k) + h(x_k)\} \) is nonincreasing,~\eqref{eq:lm:bound:step} yields
  \begin{equation}%
    \label{eq:lm:bound:step:2}
    f(x_0) + h(x_0) \geq f(x_k) + h(x_k) \geq b_f + \left( \tfrac{1}{2} \alpha (\frac{\sigma_k}{\theta_1} - \lambda^{-1}) - b_{f'} \right) \|\skcp\| + b_h.
  \end{equation}
  As \(\lim_{k \in \mathcal{K}'} \sigma_k = + \infty\), for \(k\) sufficiently large, \(\tfrac{1}{2} \alpha (\frac{\sigma_k}{\theta_1}  - \lambda^{-1}) > b_{f'}\).
  Thus, for all sufficiently large \(k \in \mathcal{K}'\),~\eqref{eq:lm:bound:step:2} combines with \(\|\skcp\| \geq \alpha\) to give
  \begin{equation*}
    f(x_0) + h(x_0) \geq b_f + \left( \tfrac{1}{2} \alpha (\frac{\sigma_k}{\theta_1}  - \lambda^{-1}) - b_{f'} \right) \alpha + b_h,
  \end{equation*}
  which is a contradiction because the right-hand side diverges.
  Thus, \(\lim_{k\in \mathcal{K}}\|\skcp\| = 0\).
  Finally, since \(\|s_k\| \leq \theta_2 \|\skcp\|\), we get also \(\lim_{k\in \mathcal{K}}\|s_k\| = 0\).
\end{proof}

For the remainder of this section, we need the following assumption.

\begin{modelassumption}%
  \label{asm:model-h-adequacy}
  For all \(k \in \N \), the model function $\psi (\cdot, x_k)$ satisfies
  \begin{equation}
    |h(x_k + s_k) - \psi(s_k; x_k)| = o(\|s_k\|) \quad \textup{as} \quad s_k \to 0.
  \end{equation}
\end{modelassumption}

\Cref{asm:model-h-adequacy} is trivially satisfied if, at each iteration $k$, we set \(\psi(s; x_k) = h(x_k + s)\), which is what \citet{kanzow-lechner-2024} do.
However, the assumption also holds when \(h(x) = g(c(x))\), where \(c : \R^n \to \R^m\) has Lipschitz-continuous or \(\alpha_h\)-H\"older-continuous Jacobian, \(g : \R^m \to \R\) is \(L\)-Lipschitz continuous, and we choose \(\psi(s; x_k) := g(c(x_k) + \nabla c(x_k)^T s)\).
Indeed, there exists \(M > 0\) such that \(|h(x_k + s) - \psi(s; x_k)| \leq L \|c(x_k + s) - c(x_k) - \nabla c(x_k)^T s\| \leq L M \|s\|^{1 + \alpha_h} = o(\|s\|)\).

\begin{theorem}%
  \label{thm:xi-bound2}
  Let \Cref{asm:g-continuous,asm:psi,asm:unif-prox-bounded,asm:model-h-adequacy} be satisfied.
  Assume that there is an index set \(\mathcal{K} \subseteq \N\) so that (i) there is \(\alpha > 0\) such that \( \nu_k ^{-1/2}\xikcp^{1/2} \ge \alpha \) for all \(k \in \mathcal{K}\), (ii) \({\{\sigma_k(1+\|B_k\|)^{-1}\}}_{k \in \mathcal{K}}\) is unbounded and (iii) \({\{x_k\}}_{k \in \mathcal{K}}\) is bounded.
  Then, there is an index set \(\mathcal{K}' \subseteq \mathcal{K}\) such that for all \(k \in \mathcal{K}' \) sufficiently large, \(k\) is a very successful iteration.
\end{theorem}

\begin{proof}
  By Assumption~(ii), there is an index set \(\mathcal{K}' \subset \mathcal{K}\) such that \(\lim_{k\in \mathcal{K}'}\sigma_k(1+\|B_k\|)^{-1} = \infty \).
  Since \(\sigma_k \ge \sigma_k(1+\|B_k\|)^{-1}\), we also have \(\lim_{k\in \mathcal{K}'}\sigma_k = \infty \).
  \Cref{lem:unbounded-sigma} then implies \(\lim_{k\in \mathcal{K}'} \|s_k\| = 0 \).
  For all \(k \in \mathcal{K}'\), \Cref{asm:model-h-adequacy} combines with~\eqref{eq:cauchy-decrease} and a Taylor expansion of \(f\) about \(x_k\) to give
  \begin{align}
    |\rho_k - 1| & = \left| \frac{ (f+h)(x_k+s_k) - (\varphi(s_k ; x_k) + \psi(s_k ; x_k)) }{ \varphi(0 ; x_k) + \psi(0 ; x_k) - (\varphi(s_k ; x_k) + \psi(s_k ; x_k)) } \right| \nonumber                                                                     \\
                 & = \left| \frac{( f+h)(x_k+s_k) - (f( x_k) + \nabla {f( x_k)}^T s_k + \tfrac{1}{2} s_k^T B_k s_k  + \psi(s_k ; x_k)) }{ \varphi(0 ; x_k) + \psi(0 ; x_k) - (\varphi(s_k ; x_k) + \psi(s_k ; x_k)) } \right|                         \nonumber \\
                 & \leq \frac{ | f(x_k+s_k) - f(x_k) - \nabla f(x_k)^T s_k |}{ (1-\theta_1) \xikcp } + \frac{ \|B_k\| \|s_k\|^2 }{ 2(1-\theta_1) \xikcp } \nonumber                                                                                             \\
                 & \phantom{\leq} \ + \frac{ | h(x_k + s_k) - \psi(s_k; x_k)| }{ (1-\theta_1) \xikcp } \nonumber                                                                                                                                                \\
                 & = \frac{ o(\|s_k\|)}{(1-\theta_1) \xikcp}  + \frac{\|B_k\| \|s_k\|^2}{2(1-\theta_1) \xikcp}  + \frac{ o(\|s_k\|) }{(1-\theta_1) \xikcp} \nonumber                                                                                            \\
                 & \leq \frac{ o(\|s_k\|)}{\xikcp}  + \frac{(1+\|B_k\|) \|s_k\|^2}{2(1-\theta_1) \xikcp}.
    \label{ineq:rho_k}
  \end{align}

  By Assumption~(i), \Cref{lem:lm-bound-Seps} implies \(\xicp(x_k; \nu_k^{-1}) \geq \frac{\alpha}{\theta_2 \sqrt{2}} \|s_k\|\) for all \(k \in \mathcal{K}'\), which we apply to the first term in the right-hand side of~\eqref{ineq:rho_k}.
  Similarly,~\eqref{eq:lm-bound-Seps2} implies
  \begin{equation*}
    \xikcp \geq \tfrac{1}{2 \theta_2^2} \nu_k^{-1} \|s_k\|^2 = \tfrac{1}{2 \theta_1 \theta_2^2} (\|B_k\| + \sigma_k ) \|s_k\|^2 \geq \tfrac{1}{2 \theta_1 \theta_2^2} \sigma_k \|s_k\|^2,
  \end{equation*}
  which we apply to the second term in the right-hand side of~\eqref{ineq:rho_k}.
  Hence,~\eqref{ineq:rho_k} simplifies to
  \begin{equation}%
    \label{eq:ineq:rho_k_simplified}
    |\rho_k - 1| \le \frac{ o(\|s_k\|) }{ \tfrac{\alpha}{\theta_2 \sqrt{2}} \|s_k\| } + \frac{ (1 + \|B_k\|) \|s_k\|^2 }{ \tfrac{(1-\theta_1)}{\theta_1 \theta_2^2}\sigma_k \|s_k\|^2} = \frac{ o(\|s_k\|) }{ \|s_k\| } + \frac{\theta_1 \theta_2^2}{(1-\theta_1)\sigma_k(1+ \|B_k\|)^{-1}}.
  \end{equation}
  By Assumption~(ii), the right-hand side of~\eqref{eq:ineq:rho_k_simplified} converges to zero.
  Thus, for all sufficiently large \(k \in \mathcal{K}'\), \( |\rho_k - 1| \le 1 - \eta_2\), which implies that \(\rho_k \geq \eta_2 \).
\end{proof}

\Cref{thm:xi-bound2} shares similarities with \citep[Theorem~\(4.1\)]{aravkin-baraldi-orban-2024} but uses weaker assumptions.
In particular, compared to \citep[Theorem~\(4.1\)]{aravkin-baraldi-orban-2024}, we do not use the Lipschitz continuity of \(\nabla f\) nor do we require model Hessians to be uniformly bounded.

\begin{lemma}%
  \label{lem:bounded-sigma}
  Let \Cref{asm:g-continuous,asm:psi,asm:unif-prox-bounded,asm:model-h-adequacy} be satisfied.
  Assume that \(\{x_k\}_{k \in \N}\) is bounded and that there is \(\alpha > 0\) such that for all \(k \in \N\), \( \nu_k^{-1/2} \xikcp^{1/2} \geq \alpha \).
  Then, \(\{\sigma_k (1+ \|B_k\|)^{-1}\}_{k\in \N}\) is bounded.
\end{lemma}

\begin{proof}
  Assume, by contradiction, that \(\{\sigma_k(1+ \|B_k\|)^{-1}\}_{k\in \N}\) is unbounded.
  Since \({\{\sigma_k\}}_{k \in \N}\) increases only on unsuccessful iterations and \(\{(1+ \|B_k\|)^{-1}\}_{k\in \N}\) is bounded, \(\{\sigma_k (1+ \|B_k\|)^{-1}\}_{k\in \mathcal{U}}\) must be unbounded, where $\mathcal{U}$ is defined in~\eqref{eq:S-U-sets}.
  Hence, using \Cref{thm:xi-bound2}, we deduce  that there is an index set \(\mathcal{U}' \subseteq \mathcal{U}\) such that for all \(k \in \mathcal{U}'\) sufficiently large, \(k\) is a very successful iteration, i.e., \(k \in \mathcal{S}\), which is absurd.
\end{proof}

Consider the following assumption

\begin{modelassumption}%
  \label{asm:unb}
  The sequence \({\{B_k\}}_{k \in \N}\) satisfies:
  \begin{equation*}
    \sum_{k \in \N } \frac{1}{r_k} = + \infty,
    \qquad
    r_k := \max_{0 \le j \le k} \|B_j\| + 1.
  \end{equation*}
\end{modelassumption}

The next theorem examines the case where \Cref{alg:R2N} generates only a finite number of successful iterations.
The proof of the second part of the theorem, which establishes stationarity of the limit point, follows ideas similar to those in \citep[Lemma~3.3]{kanzow-mehlitz-2022}.

\begin{theorem}%
  \label{thm:finite-successful}
  Let \Cref{asm:g-continuous,asm:psi,asm:unif-prox-bounded,asm:model-h-adequacy,asm:unb} be satisfied.
  If \Cref{alg:R2N} generates finitely many successful iterations, then there is \(x^* \in \R^n\) such that \(x_k = x^*\) for all sufficiently large
  \(k\), and \(\liminf_{k \to \infty} \nu_k^{-\frac{1}{2}} {\xicp(x_k, \nu_k^{-1})}^{\frac{1}{2}} = 0\).
  Moreover, \(x^*\) is stationary for~\eqref{eq:nlp}.
\end{theorem}

\begin{proof}
  First, we prove that \(\liminf_{k \to \infty} \nu_k^{-\frac{1}{2}} {\xicp(x_k, \nu_k^{-1})}^{\frac{1}{2}} = 0\).
  Assume, by contradiction, that there is \( \alpha > 0 \) such that \( \nu_k^{-\frac{1}{2}} {\xicp(x_k, \nu_k^{-1})}^{\frac{1}{2}} \geq \alpha\) for all \( k \in \N \)
  and let \(k_f\) be the last successful iteration.
  Hence, \(x_k = x_{k_f}\) for all \(k \geq k_f\) and \(\lim_k x_k = x_{k_f} = x^*\).
  Using \Cref{lem:bounded-sigma},
  \(\{\sigma_k(1+ \|B_k\|)^{-1}\}_{k \in \N}\) is bounded by a constant \(b_{\sigma}>0\).
  This, implies, that for all \(k > k_f\),
  \[
    \frac{1}{r_k} = \frac{1}{1+\max_{0 \le j \le k} \|B_j\|} \le \frac{1}{1+\|B_k\|} \le \frac{b_{\sigma}}{\sigma_k}.
  \]
  Thus, \Cref{asm:unb} implies
  \begin{equation}%
    \label{series_sigma_diverge}
    \sum_{k =k_f+1 }^{\infty} \frac{1}{\sigma_k} = + \infty.
  \end{equation}
  On the other hand, all $k > k_f$, $k$ is an unsuccessful iteration.
  The mechanism of \Cref{alg:R2N} then ensures \(\frac{\sigma_k}{\sigma_{k+1}} \le \frac{1}{\gamma_1} < 1\) for all $k > k_f$.
  But this implies that \(\sum_{k =k_f+1 }^{\infty} \frac{1}{\sigma_k}\) converges, which contradicts~\eqref{series_sigma_diverge}.
  Hence, \(\liminf_{k \to \infty} \nu_k^{-\frac{1}{2}} \xicp(x_k, \nu_k^{-1})^{\frac{1}{2}} = 0\).

  The first part of the theorem implies existence of \(\mathcal{K} \subseteq \N\) such that
  \[
    \lim_{k \in \mathcal{K}} \nu_k^{-\frac{1}{2}} \xicp(x_k, \nu_k^{-1})^{\frac{1}{2}} = 0.
  \]
  Let \(\Kf \coloneqq \{k \in \mathcal{K} \mid k > k_f\}\).
  By definition of \(\skcp\),
  \begin{equation}%
    \label{eq:def-skcp-Kf}
    \nabla f(x^*)^T \skcp + \tfrac{1}{2} \nu_k^{-1} \| \skcp \|^2 + \psi(\skcp; x^*) \leq \psi(0; x^*)
    \quad (k \in \Kf).
  \end{equation}

  Since \(\Kf \subseteq \mathcal{U}\) by assumption, \(\lim_{k \in \Kf} \sigma_k = +\infty\), and \(\{x_k\}_{k \in \Kf}\) is constant, hence bounded.
  Thus, \Cref{lem:unbounded-sigma} implies that
  \begin{equation}%
    \label{eq:skcp-converge}
    \lim_{k \in \Kf} \skcp = 0.
  \end{equation}
  Because \(\lim_{k \in \Kf} \nu_k^{-\frac{1}{2}} \xicp(x_k, \nu_k^{-1})^{\frac{1}{2}} = 0\),~\eqref{eq:xi-decrease} implies that
  \begin{equation}%
    \label{eq:nuk-skcp-converge}
    \lim_{k \in \Kf} \nu_k^{-1} \skcp = 0.
  \end{equation}

  By taking the limit superior in~\eqref{eq:def-skcp-Kf} over \(\Kf\), and using~\eqref{eq:skcp-converge} and~\eqref{eq:nuk-skcp-converge},
  \begin{equation}%
    \label{eq:psi-limsup}
    \limsup_{k \in \Kf} \psi(\skcp; x^*) \leq \psi(0; x^*).
  \end{equation}
  Using~\eqref{eq:skcp-converge} and the lower semi-continuity of \(\psi(\cdot; x^*)\) at \(0\) (from \Cref{asm:psi}), we have \(\liminf_{k \in \Kf} \psi(\skcp; x^*) \ge \psi(0; x^*)\).
  Hence,
  \begin{equation}%
    \label{eq:psi-lim}
    \lim_{k \in \Kf} \psi(\skcp; x^*) = \psi(0; x^*).
  \end{equation}

  On the other hand, since \(\skcp \in P(x_k, \nu_k^{-1})\),
  \begin{equation*}
    - \nu_k^{-1} \skcp  \in \nabla f(x^*) + \partial \psi(\skcp; x^*)
    \quad (k \in \Kf).
  \end{equation*}
  Thus, by using~\eqref{eq:nuk-skcp-converge},
  \begin{equation}%
    \label{eq:skcp-nabla-f}
    0 \in \nabla f(x^*) +\limsup_{k \in \Kf} \partial \psi(\skcp; x^*).
  \end{equation}
  Using~\eqref{eq:skcp-converge} and~\eqref{eq:psi-lim}, \citep[Proposition~8.7]{rockafellar-wets-1998} implies that
  \begin{equation*}
    \limsup_{k \in \Kf} \partial \psi(\skcp; x^*) \subseteq \partial \psi(0; x^*),
  \end{equation*}
  which, combined with~\eqref{eq:skcp-nabla-f}, gives
  \begin{equation*}
    0 \in \nabla f(x^*) + \partial \psi(0; x^*) = \nabla f(x^*) + \partial h(x^*).
  \end{equation*}
  In other words, \(x^*\) is stationary for~\eqref{eq:nlp}.
\end{proof}

To the best of our knowledge, in the case of a finite number of successful iterations, \Cref{thm:finite-successful} is the first convergence result of a regularized or trust-region method that does not rely on the boundedness of the regularization parameter or trust-region radius.
Remarkably, even in the absence of such boundedness, the theorem establishes the stationarity of the limit point \(x^*\), despite the fact that the only a subsequence of the stationarity measure \(\nu_k^{-1/2} \xicp(x_k, \nu_k^{-1})^{1/2}\) converges to zero.

Now we consider the case where the number of successful iterations is infinite.
Let $\tau \in \N_0$ and define, as in \citep{diouane-habiboullah-orban-2024},
\begin{subequations}%
  \label{eq:def-T-W}
  \begin{align}%
    \label{Tk}
    \mathcal{T}^{\tau}_k & = \left\{ j = 0, \ldots, k \mid j < \tau |\mathcal{S}_j| \right\},    \\
    \label{Wk}
    \mathcal{W}^{\tau}_k & = \left\{ j = 0, \ldots, k \mid j \geq \tau |\mathcal{S}_j| \right\}.
  \end{align}
\end{subequations}

The sets~\eqref{eq:def-T-W} are the sets defined in \citep[Equation~\((52)\)]{diouane-habiboullah-orban-2024} with \(\lambda \defeq 0\).
For readability, we drop the superscript \(\lambda\) in the following.
The next lemma provides a series comparison result that will be used in the proof of the main theorem.

\begin{lemma}[{\protect \citealp[Lemma~\(10\)]{diouane-habiboullah-orban-2024}}]%
  \label{lem:comparison-series}
  Let \({\{r_j\}}_{j \in \N}\) be a non-decreasing positive real sequence.
  For any \(k \in \N\),
  \begin{equation*}
    \tau \sum_{j \in \mathcal{S}_k} \frac{1}{r_j} \ge  \sum_{j \in \mathcal{T}^{\tau}_k} \frac{1}{r_j} = \sum_{j = 0 }^k \frac{1}{r_j} - \sum_{j \in\mathcal{W}^{\tau}_k } \frac{1}{r_j},
  \end{equation*}
  where $\mathcal{T}^{\tau}_k$ and $\mathcal{W}^{\tau}_k$ are defined in~\eqref{eq:def-T-W}.
\end{lemma}

The following lemma plays a key role in deriving a convergence result in the case where the number of successful iterations is infinite.

\begin{lemma}%
  \label{lem:series-Wk}
  Let \Cref{asm:g-continuous,asm:psi,asm:unif-prox-bounded,asm:model-h-adequacy,asm:unb} be satisfied.
  Assume that (i) \(\tau \in \N_0\) is chosen so that \(\gamma_3 \gamma_1^{\tau-1} > 1\), (ii) there is \(\alpha > 0\) such that \( \nu_k ^{-1/2}\xikcp^{1/2} \ge \alpha \) for all \(k \in \N\), and (iii) \({\{x_k\}}_{k \in \N}\) is bounded.
  Then, \(\left\{ \sum_{j \in \mathcal{W}^{\tau}_k } \frac{1}{r_j} \right\}_{k \in \N}\) is bounded, where \(r_j\) is as in \Cref{asm:unb}.
\end{lemma}

\begin{proof}
  For any \(j \geq 0\), \Cref{lem:bounded-sigma} and the update mechanism of \Cref{alg:R2N} imply that there is \(b_{\sigma} > 0\) such that
  \[
    \frac{1}{r_j} \leq \frac{1}{1 + \|B_j\|} \leq \frac{b_{\sigma}}{\sigma_j} \leq \frac{b_{\sigma}}{\gamma_3^{|S_j|} \gamma_1^{|U_j|} \sigma_0} =
    \frac{b_{\sigma}}{\gamma_3^{|S_j|} \gamma_1^{j - |S_j|} \sigma_0}.
  \]
  Consider now \(k \ge 0\) and \(j \in \mathcal{W}^{\tau}_k\).
  Then, \(j > \tau |\mathcal{S}_j|\), which, together with the fact that \(\gamma_1 > 1\) and \(0 < \gamma_3 \le 1\) leads to
  \[
    \frac{1}{r_j} \leq  \frac{b_{\sigma}}{\gamma_3^{|S_j|} \gamma_1^{j - |S_j|} \sigma_0} < \frac{b_{\sigma}}{\gamma_3^{j/\tau} \gamma_1^{j - j/\tau} \sigma_0} = \frac{b_{\sigma}}{(\gamma_3 \gamma_1^{\tau-1})^{ j / \tau} \sigma_0}.
  \]
  We sum the above inequalities over \(j \in \mathcal{W}^{\tau}_k\) and use the fact that \(\gamma_3 \gamma_1^{\tau-1} > 1\) to obtain
  \begin{align*}%
    \sum_{j \in \mathcal{W}^{\tau}_k } \frac{1}{r_j} <  \frac{b_{\sigma}}{\sigma_0}\sum_{j \in \mathcal{W}^{\tau}_k } \frac{1}{(\gamma_3 \gamma_1^{\tau-1})^{ j / \tau}} \le \frac{b_{\sigma}}{\sigma_0}\sum_{j \in \N} \frac{1}{(\gamma_3 \gamma_1^{\tau-1})^{ j / \tau}} < \infty.
  \end{align*}
\end{proof}

We state now our main convergence result.

\begin{theorem}%
  \label{thm:R2N-unbounded}
  Let \Cref{asm:g-continuous,asm:psi,asm:unif-prox-bounded,asm:model-h-adequacy,asm:unb} be satisfied.
  Assume that \Cref{alg:R2N} generates infinitely many successful iterations and that there is \({(f+h)}_{\text {low }} \in \mathbb{R}\) such that \((f+h)\left(x_k\right) \ge{(f+h)}_{\mathrm{low }}\) for all \(k \in \mathbb{N}\).
  Then, \(\liminf_{k \to +\infty} \nu_k ^{-1/2} \xikcp^{1/2} = 0\).
\end{theorem}

\begin{proof}
  By contradiction, assume that there is \(\alpha>0\) such that for all \(k \in \N\), \(\nu_k^{-1/2} {\xikcp}^{1/2} \geq \alpha\).
  \Cref{lem:xk-convergennt} shows that \({\{x_k\}}_{k \in \N}\) is convergent, hence bounded.
  \Cref{lem:bounded-sigma} then implies that \({\{\sigma_k(1+ \|B_k\|)^{-1}\}}_{k \in \N}\) is bounded, say by \(b_{\sigma} > 0\).
  Equivalently, \(\sigma_k \le b_{\sigma} (1+ \|B_k\|)\).
  For any \(j \in \mathcal{S}\), \(\rho_j \ge \eta_1\), and
  \begin{align*}
    f(x_j) + h(x_j) - f(x_{j+1}) - h(x_{j+1}) & \geq \eta_1 (1 - \theta_1) \xicp(x_j, \nu_j^{-1})
    \\
                                              & \geq \eta_1 (1 - \theta_1) \nu_j \alpha^2
    \\
                                              & = \frac{\eta_1 (1 - \theta_1) \theta_1 \alpha^2}{\sigma_j + \|B_j\|}
    \\
                                              & \geq \frac{\eta_1 (1-\theta_1)\theta_1 \alpha^2}{(b_{\sigma} (1 + \|B_j\|) + \|B_j\|)}
    \\
                                              & \geq \frac{\eta_1 (1-\theta_1)\theta_1 \alpha^2}{(1 + b_{\sigma}) (1 + \|B_j\|)}
    \\
                                              & \geq \frac{\eta_1 (1-\theta_1)\theta_1 \alpha^2}{1 + b_{\sigma}} \ \frac{1}{r_j},
  \end{align*}
  where \(r_j\) is defined in \Cref{asm:unb}.
  Let \(k \in \mathcal{S}\).
  We sum the above inequalities over all \(j \in \mathcal{S}_k\), and obtain
  \begin{equation*}
    f(x_0) + h(x_0) - f(x_{k+1}) - h(x_{k+1}) 
    \geq \frac{\eta_1 (1-\theta_1) \theta_1 \alpha^2}{1 + b_\sigma} \sum_{j \in \mathcal{S}_k} \frac{1}{r_j}.
  \end{equation*}
  Since \(f + h\) is bounded below, it follows that \(\sum_{k \in \mathcal{S}} \frac{1}{r_k} < \infty\).
  Let \(\tau \in \N_0\) be chosen so that \(\gamma_3 \gamma_1^{\tau-1} > 1\).
  By \Cref{lem:series-Wk}, \(\sum_{j \in \mathcal{W}^{\tau}_k } \frac{1}{r_j}\) is uniformly bounded for all \(k \in \N\).
  However, \Cref{lem:comparison-series} yields that for all \(k \in \N\),
  \begin{equation*}
    \sum_{j = 0}^k \frac{1}{r_j} \leq \tau \sum_{j \in \mathcal{S}_k} \frac{1}{r_j} + \sum_{j \in \mathcal{W}^{\tau}_k } \frac{1}{r_j},
  \end{equation*}
  which implies that \(\sum_{k = 0}^{\infty} \frac{1}{r_k}\) converges, and contradicts \Cref{asm:unb}.
\end{proof}

Note that the assumptions involved in \Cref{thm:R2N-unbounded} are weaker compared to existing methods in the literature---see \Cref{sec:introduction}.

Moreover, our result is more general in that it establishes the existence of a subsequence \({(x_k)}_{k \in \mathcal{K}}\) such that \(\lim_{k \in \mathcal{K}} \nu_k^{-1/2} \, \xi_{\textup{cp}}(\skcp; x_k, \nu_k^{-1})^{1/2} = \lim_{k \in \mathcal{K}} \nu_k^{-1} \, \|\skcp\| = 0\).
That is in contrast with \citet{kanzow-mehlitz-2022}, who prove that \(\lim_{k \in \mathcal{K}} \nu_k^{-1} \, \|\skcp\| = 0\) only for subsequences \(\mathcal{K}\) along which \({(x_k)}_{k \in \mathcal{K}}\) converges.
Furthermore, their analysis relies on the assumption that \((\nu_k^{-1})_{k \in \N}\) is bounded away from zero, which is not required in our case.

Regarding the stationarity of accumulation points of \Cref{alg:R2N}, \citep[Theorem~5]{leconte-orban-2023-2} directly implies that if there exists a subsequence \(\{x_k\}_{k \in \mathcal{K}} \to x^*\) such that \(\lim_{k \in \mathcal{K}} \nu_k^{-1/2} \, \xi_{\textup{cp}}(\skcp; x_k, \nu_k^{-1})^{1/2} = 0\), and such that there exists a limiting model \(\psi(\cdot; x^*)\) that satisfies \citep[Model Assumption~3.1]{leconte-orban-2023-2}, then \(x^*\) is a stationary point of~\eqref{eq:nlp}, provided that the sequence \((\sigma_k)_{k \in \N}\) is bounded from below.

The analysis of \citep{leconte-orban-2023-2} does not establish whether all accumulation points of \(\{x_k\}_{k \in \N}\) are stationary for~\eqref{eq:nlp}.
A complete characterization of the stationarity of all accumulation points of \Cref{alg:R2N} under similar or weaker assumptions, in particular the removal of the boundedness assumption on the sequence \((\sigma_k)_{k \in \N}\), is left for future work.

\section{Complexity analysis of \texorpdfstring{\Cref{alg:R2N}}{Algorithm 4.1}}%
\label{sec:evaluation-complexity}

In this section, we study the evaluation complexity of \Cref{alg:R2N} in the case where the model Hessians are allowed to be unbounded.
We replace \Cref{asm:model-h-adequacy} with the following.

\begin{modelassumption}%
  \label{asm:model-adequacy2}
  There is \(\kappa_{\textup{m}} > 0\) such that for all \(k \in \N\),
  \begin{equation}
    |(f + h)(x_k + s_k) - (\varphi +
    \psi)(s_k; x_k)| \leq \kappa_{\textup{m}} (1 + \|B_k\|) \|s_k\|^2.
  \end{equation}
\end{modelassumption}

Note that if \(\nabla f\) is Lipschitz continuous and \(\psi(;x)\) satisfies \Cref{asm:model-h-adequacy} then \Cref{asm:model-adequacy2} is satisfied as discussed by \citet{leconte-orban-2023-2}.

The next lemma will allow us to show that \(\{\sigma_k (1+\|B_k\|^{-1})\}\) is bounded

\begin{lemma}%
  \label{thm:success-bound}
  Let \Cref{asm:g-continuous,asm:psi,asm:unif-prox-bounded,asm:model-adequacy2} be satisfied.
  Define
  \begin{equation*}
    b_{\textup{succ}} := \frac{2\kappa_{\textup{m}}}{1 - \eta_2} > 0.
  \end{equation*}
  If \(x_k\) is not first-order stationary and \(\sigma_k (1 + \max_{0 \le j \le k} \|B_j\|)^{-1} \geq b_{\textup{succ}}\), iteration \(k\) is very successful and \(\sigma_{k+1} \le \sigma_k\).
\end{lemma}

\begin{proof}
  Because \(x_k\) is not first-order stationary, \(s_k \neq 0\).
  By definition of \(s_k\), \( m(0; x_k,\sigma_k) \ge m(s_k;x_k,\sigma_k)\).
  Hence,
  \begin{equation}%
    \label{eq:bound_xicp}
    \varphi_k\left(0 ; x_k\right) +\psi_k\left(0 ; x_k\right) \geq \varphi_k\left(s_k ; x_k\right) +\psi_k\left(s_k ; x_k\right) + \tfrac{1}{2}\sigma_k \|s_k\|^2.
  \end{equation}
  Thus, \Cref{asm:model-adequacy2} yields
  \begin{align*}
    |\rho_k - 1| & = \left|\frac{( f+h)(x_k+s_k) - (\varphi_k(s_k ; x_k) + \psi_k(s_k ; x_k)) }{ \varphi_k(0 ; x_k) + \psi_k(0 ; x_k) - (\varphi_k(s_k ; x_k) + \psi_k(s_k ; x_k)) }\right| \\
                 & \le \frac{ \kappa_{\textup{m}} (1 + \|B_k\|) \|s_k\|^2 }{ \tfrac{1}{2} \sigma_k \|s_k\|^2 }                                                                              \\
                 & \le \frac{ 2\kappa_{\textup{m}} }{ \sigma_k (1 + \max_{0 \le j \le k} \|B_j\|)^{-1}}
    \le \frac{ 2\kappa_{\textup{m}} }{ b_{\textup{succ}} } = 1-\eta_2.
  \end{align*}
  Thus, we obtain \(\rho_k \ge \eta_2\), meaning that the iteration $k$ is very successful.
\end{proof}

The next theorem shows that \(\{\sigma_k (1+\|B_k\|^{-1})\}\) is bounded.

\begin{theorem}%
  \label{thm:sigma-max}
  Let \Cref{asm:g-continuous,asm:psi,asm:unif-prox-bounded,asm:model-adequacy2} be satisfied.
  For all \(k \in \N\), if \(x_k\) is not stationary,
  \begin{equation*}
    \sigma_k (1 + \max_{0 \le j \le k} \|B_j\|)^{-1} \le b_{\max} := \min \, \{  \sigma_0 (1+ \|B_0\|)^{-1}, \, \gamma_2 b_{\textup{succ}} \} > 0.
  \end{equation*}
\end{theorem}

\begin{proof}
  Set \(b_k := \sigma_k (1 + \max_{0 \le j \le k} \|B_j\|)^{-1}\) for all \(k\).
  We proceed by induction.
  For \(k = 0\), \( \sigma_0 (1+ \|B_0\|)^{-1} \le b_{\max}\) by definition.
  Assume that \(b_k \le b_{\max}\) for \(k \geq 0\).

  Assume first that \(b_k < b_{\textup{succ}}\).
  Because \(\{(1 + \max_{0 \le j \le k} \|B_j\|)^{-1}\}\) is non-increasing, the update of \(\sigma_k\) in \Cref{alg:R2N} ensures that
  \begin{equation*}
    b_{k+1} = (1 + \max_{0 \le j \le k+1} \|B_j\|)^{-1} \sigma_{k+1} \le (1 + \max_{0 \le j \le k} \|B_j\|)^{-1} \gamma_2 \sigma_k = \gamma_2 b_{k} < \gamma_2 b_{\textup{succ}} \le b_{\max}.
  \end{equation*}

  Now, assume conversely that \(b_k \ge  b_{\textup{succ}}\).
  \Cref{thm:success-bound} implies that iteration \(k\) is very successful, and \(\sigma_{k+1} \le \sigma_k\).
  Thus,
  \begin{align*}
    b_{k+1} = (1 + \max_{0 \le j \le k+1} \|B_j\|)^{-1} \sigma_{k+1} < (1 + \max_{0 \le j \le k} \|B_j\|)^{-1}  \sigma_{k} = b_k \le b_{\max}.
    \tag*{\qed}
  \end{align*}
\end{proof}

Additionally, instead of \Cref{asm:unb}, we assume that model Hessians grow at most linearly with \(|\mathcal{S}_k|\), which covers multiple quasi-Newton approximations---see \Cref{sec:introduction}.

\begin{modelassumption}%
  \label{asm:unb-success}
  There are \(\mu > 0\) and \(0 \leq p \leq 1\) such that, for all \(k \in \N\),
  \begin{equation}%
    \label{eq1:qsm:unb}
    \max_{0 \le j \le k} \|B_j\| \le \mu (1 + |\mathcal{S}_k|^p).
  \end{equation}
\end{modelassumption}
Because \(|\mathcal{S}_k|\) is non-decreasing with \(k\),~\eqref{eq1:qsm:unb} is equivalent to \(\|B_k\| \le \mu(1 + |\mathcal{S}_k|^p)\) for all  \(k \in \N\).
The following theorem considers the case with a finite number of successful iterations.
The proof follows \citep[Theorem~\(4.2\)]{aravkin-baraldi-orban-2024} and is recalled here for completeness.

\begin{theorem}%
  \label{thm:finite-successful-Lipschitz}
  Let \Cref{asm:g-continuous,asm:psi,asm:unif-prox-bounded,asm:model-adequacy2,asm:unb-success} be satisfied.
  If \Cref{alg:R2N} generates finitely many successful iterations, then \(x_k = x^*\) for all sufficiently large \(k\) where \(x^*\) is a stationary point.
\end{theorem}

\begin{proof}
  Assume by contradiction that \(x^*\) is not a stationary point.
  Because the number of successful iterations is finite, according to \Cref{asm:unb-success}, there is \(k_f \in \N\) such that \(\|B_k\| \leq \mu (1 + |\mathcal{S}_{k_f}|^p)\) for all \(k \geq k_f\), where \(k_f\) is the index of the last successful iteration.
  The mechanism of \Cref{alg:R2N} ensures that \(\sigma_k\) increases on unsuccessful iterations.
  Hence, there must exist an unsuccessful iteration \(k > k_f\) such that \(\sigma_{k} \ge b_{\textup{succ}} (1 + \mu (1 + |\mathcal{S}_{k_f}|^p)) \ge b_{\textup{succ}} (1 + \|B_{k}\|)\), with \(b_{\textup{succ}}\) defined in \Cref{thm:success-bound}.
  Because \(x^*\) is not stationary, we can apply \Cref{thm:success-bound}, which shows that \(k\) is very successful, and contradicts our assumption.
\end{proof}

We know from \Cref{thm:R2N-unbounded} that \(\liminf_{k \to +\infty} \nu_k ^{-1/2} \xikcp^{1/2} = 0\) when \Cref{alg:R2N} generates infinitely many successful iterations.
Let \(\epsilon > 0\) and \(k_{\epsilon}\) be the first iteration of \Cref{alg:R2N} such that \(\nu_k ^{-1/2}\xikcp^{1/2} \leq \epsilon\).
Define
\begin{subequations}%
  \label{eq:eps:S-U-sets-unbounded}
  \begin{align}
    \mathcal{S}(\epsilon) & := \mathcal{S}_{k_\epsilon-1} = \{ k \in \mathcal{S} \mid k < k_{\epsilon} \},
    \\
    \mathcal{U}(\epsilon) & := \mathcal{U}_{k_\epsilon-1} = \{k \in \N \mid k \not \in \mathcal{S} \text{ and } k < k_{\epsilon} \}.
  \end{align}
\end{subequations}

The next theorems bound \(k_\epsilon\).
The proofs are similar to \citep[Theorem 2]{diouane-habiboullah-orban-2024}.

\begin{theorem}%
  \label{thm:complexity:S}
  Let \Cref{asm:g-continuous,asm:psi,asm:unif-prox-bounded,asm:model-adequacy2,asm:unb-success} be satisfied.
  Assume that \Cref{alg:R2N} generates infinitely many successful iterations and that there is \({(f + h)}_{\mathrm{low }} \in \mathbb{R}\) such that \((f + h)\left(x_k\right) \ge {(f + h)}_{\mathrm{low }}\) for all \(k \in \mathbb{N}\).
  If \(0 \le p < 1\),
  \begin{equation}%
    \label{eq:S-eps2}
    |\mathcal{S}(\epsilon)| \leq {((1 - p) \kappa_1 \epsilon^{-2} + 1)}^{1 /(1- p)} - 1 = O(\epsilon^{-2 /(1- p)}),
  \end{equation}
  where
  \begin{equation*}
    \kappa_1 = \frac{ ((f+h)\left(x_0\right)-{(f+h)}_{\mathrm{low }}) \left(b_{\max} +2 \mu (1+b_{\max})\right)} {\eta_1 \theta_1 (1-\theta_1)},
  \end{equation*}
  and \(b_{\max}\) is as in \Cref{thm:sigma-max}.
  If \(p = 1\),
  \begin{equation}%
    \label{eq:S-exp}
    |\mathcal{S}(\epsilon)| \leq \exp(\kappa_1\epsilon^{-2}) - 1.
  \end{equation}
\end{theorem}

\begin{proof}
  Let \(k \in \mathcal{S}(\epsilon)\), then \(\nu_k ^{-1/2}\xikcp^{1/2} \ge  \epsilon\) and
  \begin{equation}%
    \label{eq:cauchydec}%
    (f + h)
    \left(x_k\right)-(f + h)\left(x_k + s_k\right)  \ge \eta_1 (1-\theta_1) \xicp(x_k; \nu_k^{-1})  \ge \eta_1 (1-\theta_1) \nu_k \epsilon^2.
  \end{equation}
  \Cref{thm:sigma-max} implies
  \begin{align*}
    \nu_k = \frac{\theta_1}{\|B_k\| + \sigma_k} & \ge \frac{\theta_1}{\max_{0 \le j \le k} \|B_j\| + b_{\max} (1 + \max_{0 \le j \le k} \|B_j\|)}
    \\
                                                & = \frac{\theta_1}{b_{\max} +  (1 +b_{\max}) \max_{0 \le j \le k} \|B_j\|}.
  \end{align*}
  \Cref{asm:unb-success} then implies
  \begin{equation}%
    \label{eq:bound:nuk}
    \nu_k \ge \frac{\theta_1}{b_{\max} + \mu (1 +b_{\max})(1+|\mathcal{S}_k|^p)} = \frac{\theta_1}{|\mathcal{S}_k|^p} \zeta\left(|\mathcal{S}_k|^p\right),
  \end{equation}
  where \(\zeta: \R_+ \to \R\), \(\zeta(x) := x / (b_{\max} + \mu (1 + b_{\max}) (x + 1))\).

  Because \(\zeta\) is non-decreasing and \(|\mathcal{S}_k|\ge 1\) (as we have infinitely many successful iterations), \(\zeta(|\mathcal{S}_k|^p) \ge \zeta(1) = \left(1 +2 \mu (1+b_{\max}) \right)^{-1} \).
  Thus,~\eqref{eq:bound:nuk} becomes
  \begin{equation*}
    \nu_k \ge   \frac{\theta_1}{b_{\max} +2 \mu (1+b_{\max})} \frac{1}{|\mathcal{S}_k|^p},
  \end{equation*}
  which combines with~\eqref{eq:cauchydec} to yield
  \begin{equation}%
    \label{eq:cauchydec:Sk}%
    (f + h)
    \left(x_k\right)-(f + h)\left(x_k + s_k\right)  \ge \frac{\eta_1 \theta_1 (1-\theta_1)\epsilon^2 }{b_{\max} +2 \mu (1+b_{\max})} \frac{1}{|\mathcal{S}_k|^p} := C \frac{1}{|\mathcal{S}_k|^p}.
  \end{equation}
  We sum over all \(k \in \mathcal{S}(\epsilon)\), and obtain
  \begin{equation*}
    (f+h)(x_0) - {(f+h)}_{\textup{low}} \ge C \sum_{k \in \mathcal{S}(\epsilon)} \frac{1}{|\mathcal{S}_k|^p} = C \sum_{k = 0}^{|\mathcal{S}(\epsilon)| - 1}\frac{1}{|\mathcal{S}_{\phi(k)}|^{p}},
  \end{equation*}
  where \(\phi\) is an increasing map from \(\{0, \ldots, |\mathcal{S}(\epsilon)| - 1\}\) to \(\mathcal{S}(\epsilon)\).
  Thus, by definition of \(\phi\) and \(\mathcal{S}_{\phi(k)}\), \(|\mathcal{S}_{\phi(k + 1)}| = |\mathcal{S}_{\phi(k)}| + 1\) and \(|\mathcal{S}_{\phi(0)}| = 1\).
  In other words, \(|\mathcal{S}_{\phi(k)}| = k + 1\), and
  \begin{equation*}
    (f+h)(x_0) - {(f+h)}_{\textup{low}} \ge C \sum_{k = 0}^{|\mathcal{S}(\epsilon)| - 1} \frac{1}{{(k + 1)}^{p}} = C \sum_{k = 1}^{|\mathcal{S}(\epsilon)|} \frac{1}{k^{p}}.
  \end{equation*}
  Because \(\int_k^{k + 1} \frac{1}{t^p} {\mathrm d}t \le \int_k^{k + 1} \frac{1}{k^p} {\mathrm d}t = \frac{1}{k^p}\),
  \begin{equation}
    (f+h)(x_0) - {(f+h)}_{\textup{low}} \ge C \sum_{k = 1}^{|\mathcal{S}(\epsilon)|} \int_k^{k + 1} \frac{1}{t^p} {\mathrm d}t = C \int_{1}^{|\mathcal{S}(\epsilon)| + 1} \frac{1}{t^p} {\mathrm d}t.
  \end{equation}
  There are two cases to consider:

  \begin{itemize}
    \item if \(0 \le p < 1\),
      \(
      (f+h)(x_0) - {(f+h)}_{\textup{low}} \geq C \frac{{(|\mathcal{S}(\epsilon)| + 1)}^{1-p} - 1}{1-p}
      \),
      which is~\eqref{eq:S-eps2};
    \item if \(p = 1\),
      \(
      (f+h)(x_0) - {(f+h)}_{\textup{low}} \geq C \log(|\mathcal{S}(\epsilon)| + 1)
      \),
      which is~\eqref{eq:S-exp}.
  \end{itemize}
\end{proof}

Finally, we derive a bound on the cardinality of \(\mathcal{U}(\epsilon)\).

\begin{theorem}%
  \label{thm:complexity:U}
  Let \Cref{asm:g-continuous,asm:psi,asm:unif-prox-bounded,asm:model-adequacy2} hold.
  Assume that \Cref{alg:R2N} generates infinitely many successful iterations.
  Then
  \begin{align}%
    \label{eq:complexity-3}
    |\mathcal{U}(\epsilon)| & \leq |\log_{\gamma_1} (\gamma_3)| |\mathcal{S}(\epsilon)| + \log_{\gamma_1}(1 + \mu(1 + |\mathcal{S}(\epsilon)|^p)) + \frac{\log(b_{\max} / \sigma_0)}{\log(\gamma_1)},
  \end{align}
  where \(\mu\) and \(p\) are defined in \Cref{asm:unb}, \(b_{\max}\) as in \Cref{thm:sigma-max}, and \(|\mathcal{S}(\epsilon)|\) is as in \Cref{thm:complexity:S}.
\end{theorem}
\begin{proof}
  The mechanism of \Cref{alg:R2N} guarantees that for all \(k \in \N\),
  \(
  |\mathcal{U}_k| \leq |\log_{\gamma_1} (\gamma_3)| \, |\mathcal{S}_k| +  \log_{\gamma_1}(\sigma_k / \sigma_0)
  \).
  Hence, \Cref{thm:sigma-max} yields
  \begin{align*}
    |\mathcal{U}(\epsilon)| & \leq |\log_{\gamma_1} (\gamma_3)| |\mathcal{S}(\epsilon)| + \log_{\gamma_1}\left(\frac{b_{\max} (1 + \max_{0 \le j \le k_\epsilon -1}\|B_j\|)}{\sigma_0}\right)              \\
                            & \leq |\log_{\gamma_1} (\gamma_3)| |\mathcal{S}(\epsilon)| + \log_{\gamma_1}(1 + \mu(1 + |\mathcal{S}(\epsilon)|^p)) + \log_{\gamma_1}\left(\frac{b_{\max}}{\sigma_0}\right).
    \tag*{\qed}
  \end{align*}
\end{proof}

The complexity bound in \Cref{thm:complexity:S} is of the same order as that of \citep[Lemma~\(4.3\)]{aravkin-baraldi-orban-2024} for trust-region methods when \(p = 0\) in \Cref{asm:unb-success}, which corresponds to bounded model Hessians.
Unlike \citep[Lemma~\(3.6\)]{aravkin-baraldi-orban-2022}, the constant \(\theta_2\), as defined in the switch on \Cref{alg:R2N:step-comparison} of \Cref{alg:R2N}, does not appear in our complexity bound.
Thus, large values of \(\theta_2\) in \Cref{alg:R2N} will not worsen the complexity bound.
In the general case where \(p > 0\), our bound is better than that in \citep[Theorem~\(4.2\)]{leconte-orban-2023-2}, as their step computation rule makes the bound dependent on \(\theta_2\).
As \(p\) approaches \(1\), the bound in \citep[Theorem~\(4.2\)]{leconte-orban-2023-2} goes to infinity, whereas ours, though exponential, remains finite, as in \citep{diouane-habiboullah-orban-2024}.
Finally, the same example as in \citep[\S \(3.1\)]{diouane-habiboullah-orban-2024} shows that our complexity bounds are also tight.

\section{Algorithmic refinements}%
\label{sec:refinements}

We describe a special case of \Cref{alg:R2N} and an extension for which the convergence theory continues to hold, that we exploit in the numerical experiments of \Cref{sec:numerical}, and that prove to be efficient in practice.
As both refinements have already been studied by \citet{leconte-orban-2025} in the context of their trust-region method, we keep our description to a minimum.

\subsection{Special case: diagonal model Hessians}%
\label{sec:dh}

If we select \(B_k\) to be diagonal in \Cref{alg:R2N}, a specialized implementation emerges whenever \(h\) is separable and \(\psi(\cdot; x_k)\) is chosen to be separable at each iteration.
For a number of choices of separable \(h\) that are of interest in applications, the step \(s_k\) may be computed analytically without requiring an iterative subproblem solver.
We refer to this implementation as R2DH, where ``DH'' stands to \emph{diagonal Hessians}.
This section is modeled after \citep[Section~\(4\)]{leconte-orban-2025}, to which we refer the reader for further information.

Diagonal quasi-Newton methods originate from \citep{gilbert-lemarechal-1989,dennis-wolkowicz-1993,nazareth-1995}.
In order for a variational problem to possess a solution that defines a diagonal update, the classic secant equation is replaced with the \emph{weak} secant equation \(s_k^T B_{k+1} s_k = s_k^T y_k\), where \(y_k = \nabla f(x_{k+1}) - \nabla f(x_k)\).
A handful of diagonal updates have been proposed in the literature.
The most efficient is probably the \emph{spectral} update \(B_{k+1} = \tau_{k+1} I\), where \(\tau_{k+1} := s_k^T y_k / s_k^T s_k\) is defined as in the spectral gradient method \citep{birgin-martinez-raydan-2014}.
Because \(B_k\) is a multiple of the identity, \(h\) and \(\psi(\cdot; x_k)\) need not be separable as the computation of \(s_k\) boils down to the evaluation of a proximal operator with step length \(1 / \sqrt{\tau_k + \sigma_k}\)---see \Cref{def:prox}.
\citet{zhu-nazareth-wolkowicz-1999} derive an update akin to the well-known PSB formula that may be indefinite.
We refer to it below as \emph{PSB}.
\citet{andrei-2019} derives an update based on a different variational problem that may also be indefinite.
We refer to it below as \emph{Andrei}.
Additionally, we include a new diagonal variant inspired from the BFGS formula using a diagonal update.
The main idea comes from applying \citep[Lemma~\(5.1\)]{chouzenoux-pesquet-repetti-2014} to the last term of the BFGS update, i.e., $y_k y_k^T / s_k^T y_k$, to obtain the diagonal update
\[
  D_{k+1} = \frac{ \sum_{i=1}^n \left| (y_k)_i\right|}{s_k^{T} y_k} \operatorname{diag}(| y_k |).
\]
This update remains positive as long as \(s_k^{T} y_k > 0\).
We refer to this variant below as \emph{DBFGS}.
Although DBFGS does not always satisfy the secant equation, our numerical results demonstrate its competitiveness against other state-of-the-art diagonal-based methods.
Note that the three updates (i.e., PSB, Andrei and DBFGS) generate \(B_k\) that is not a multiple of the identity, and hence \(h\) and \(\psi(\cdot; x_k)\) should be separable.

R2DH may act as standalone solver for~\eqref{eq:nlp} or as subproblem solver to compute \(s_k\) in \Cref{alg:R2N}.
Our results in \Cref{sec:numerical} illustrate that, in both use cases, R2DH typically outperforms R2 \citep[Algorithm~\(6.1\)]{aravkin-baraldi-orban-2022,aravkin-baraldi-leconte-orban-2021}.

\subsection{Non-monotone variants}%
\label{sec:nonmonotone}

Inspired by the success of the non-monotone spectral gradient method \citep{birgin-martinez-raydan-2014}, \citet[Section~\(6\)]{leconte-orban-2025} explain how to modify an algorithm similar to \Cref{alg:R2N} to incorporate a non-monotone strategy.

Let \(q \in \N\) be a given \emph{memory parameter}.
Define \(q_k = 1\) if \(q = 0\) and \(q_k := \min(k, \, q)\) if \(q > 0\),
Define also \(\mathcal{S}_{q_k}^+\) the set of the \(q_k\) \emph{most recent} successful iterations.
By convention, we set \(\mathcal{S}_0^+ = \{0\}\).
An iteration \(k\) now considers the objective value at each iteration in \(\mathcal{S}_{q_k}^+\).
Define
\begin{equation}%
  \label{eq:obj-max}
  (f + h)_{\max, k}%
  := \max\{ (f + h)(x_j) \mid j \in \mathcal{S}_{q_k}^+ \}.
\end{equation}
\Cref{alg:R2N} corresponds to \(q = 0\).
The non-monotone strategy consists in enforcing decrease with respect to \((f + h)_{\max, k}\) instead of \((f + h)(x_k)\).
In other words, we redefine
\[
  \rho_k :=
  \frac{
    (f + h)_{\max, k} - (f + h)(x_k + s_k)
  }{
    (f + h)_{\max, k} - (\varphi + \psi)(s_k; x_k)
  }.
\]
As in \citep[Section~\(6\)]{leconte-orban-2025}, the new expression of \(\rho_k\) does not interfere with convergence properties or complexity bounds, except that it changes the constants in the latter.
This is a positive result, especially in comparison to \citep{kanzow-mehlitz-2022}, where additional assumptions, including uniform continuity of \(f + h\), are required to establish convergence of their non-monotone proximal-gradient method.

\section{Numerical experiments}%
\label{sec:numerical}

Our implementation of \Cref{alg:R2N} and all solvers used in the experiments are available in the \href{https://github.com/JuliaSmoothOptimizers/RegularizedOptimization.jl}{RegularizedOptimization} Julia module \citep{baraldi-leconte-orban-regularized-optimization-2024,diouane-gollier-habiboullah-orban-2025}.
We compare the performance of R2N and its variants against PANOC \citep{stella-themelis-sopasakis-patrinos-2017}.
By default, R2N uses an L-BFGS approximation with memory \(5\), as implemented in the \href{https://github.com/JuliaSmoothOptimizers/LinearOperators.jl}{LinearOperators} Julia module \citep{Leconte_LinearOperators_jl_Linear_Operators_2023}, and uses parameters \(\theta_1 = (1 + \varepsilon_M^{1/5})^{-1} \approx 0.999\), \(\theta_2 = 1/\varepsilon_M \approx 10^{15}\), \(\eta_1 = \varepsilon_M^{1/4} \approx 10^{-4}\), \(\eta_2 = 0.9\), and \(\sigma_0 = \varepsilon_M^{1/3} \approx 10^{-6}\), where \(\varepsilon_M\) is the machine epsilon.
The reason for defining values based on \(\varepsilon_M\) is that our code may be run in various floating-point arithmetics.
Here, however, all tests are run in double precision.
If iteration \(k\) of \Cref{alg:R2N} is very successful,  \(\sigma_{k+1} = \sigma_k / 3\); if iteration \(k\) is unsuccessful, \(\sigma_{k+1} = 3 \sigma_k\).
Otherwise, \(\sigma_{k+1} = \sigma_k\).

All solvers use the same stopping criterion and terminate when
\begin{equation}%
  \label{eq:r2n-stop-test}
  \nu_k^{-1} \|s_{k,\textup{cp}}\| < \epsilon_a,
\end{equation}
where \(\epsilon_a = \varepsilon_M^{3/10} \approx 2 \cdot 10^{-5}\) is an absolute tolerance, or exceed the budget of \(5,000\) iterations and \(3,600\) seconds of CPU time.
This absolute stopping criterion~\eqref{eq:r2n-stop-test} aligns with the one used by default in PANOC \citep{stella-themelis-sopasakis-patrinos-2017}.
As discussed in \Cref{sec:models}, \(\nu_k^{-1} \|s_{k,\textup{cp}}\|\) can also be considered as a first-order stationarity measure for~\eqref{eq:nlp}.
\Cref{thm:R2N-unbounded} and\eqref{eq:xi-norm} imply that \(\liminf_{k \to +\infty} \nu_k^{-1} \|s_{k,\textup{cp}}\| = 0\).
To solve the subproblem in Line~\ref{alg:R2N:step-computation} of \Cref{alg:R2N}, we use either R2 \citep[Algorithm~\(6.1\)]{aravkin-baraldi-orban-2022}, or one of several R2DH variants (Spec, PSB, Andrei, or DBFGS) as described in \Cref{sec:refinements}, as well as the non-monotone  spectral R2DH (R2DH-Spec-NM) with memory \(5\).
R2 initializes \(\nu_0 = 1.0\).
R2N and R2DH initialize \(\nu_0\) according to Line~\ref{alg:R2N:step-nuk} of \Cref{alg:R2N}.
The subproblem solvers terminate as soon as
\[
  \hat{\nu}_k^{-1/2} \hat{\xi}_{\mathrm{cp}}\left( x_k+s, \hat{\nu_k} \right)^{1/2} \leq
  \begin{cases}
    10^{-3}                                                                                               & \text{if } k = 0, \\
    \min \left(\left(\nu_k^{-1} \xicp\right)^{3/4}, \, 10^{-3} \left(\nu_k^{-1} \xicp\right)^{1/2}\right) & \text{if } k > 0,
  \end{cases}
\]
where \(\xicp =  \xicp\left( x_k, \nu_k \right)\), \(\hat{\nu}_k\) and \(\hat{\xi}_{\mathrm{cp}}\) are the step size and first-order stationarity measure related to the subproblem solver.
Note that R2 and all the R2DH variants can also be used to solve~\eqref{eq:nlp} directly.
PANOC \citep{stella-themelis-sopasakis-patrinos-2017} is run with all default parameters.
All quasi-Newton approximations are initialized to the identity.
In all experiments, we use \(\psi (s; x) := h(x + s)\).

Our objective is to minimize the number of objective and gradient evaluations, as they are generally expensive to compute, while assuming that the proximal operators of common regularizers such that \(\ell_0 \) and \(\ell_1\) norms are comparatively cheap to evaluate.
We include also other test problems with the nuclear norm and the rank regularizers.

In our figures, we set \((f+h)^*\) to the best value found by all the solvers.
We plot \(\Delta (f+h) (x_k)=(f+h)(x_k) - (f+h)^*\) against the iterations to illustrate progress towards that best value.
We also report the following solver statistics in tables: the final value of \(f\) at convergence; the final \(h / \lambda\), where \(\lambda\) is a weight on the regularizer \(h\); the final stationarity measure \(\nu^{-1}\|s\|\); the number of evaluations of the smooth objective (\(\# f\)); the number of evaluations of the gradient (\(\# \nabla f\)); the number of proximal operator evaluations (\(\# \text{prox}\)); and the elapsed time \(t\) in seconds.

\subsection{Basis pursuit denoise (BPDN)}

The first set of experiments focuses on the basis pursuit denoise problem as described in \citep{aravkin-baraldi-orban-2022}, which is common in statistical and compressed sensing applications.
The goal is to recover a sparse signal \(x_{\text{true}} \in \mathbb{R}^n\) from noisy observed data \(b \in \mathbb{R}^m\).
This problem can be formulated as
\begin{equation}%
  \label{eq:bpdn}
  \minimize{x} \ \tfrac{1}{2} \|A x - b\|_2^2 + \lambda \|x\|_0,
\end{equation}
where \(A\) is \(m\)\(\times\)\(n\) and randomly generated with orthonormal rows.
We set \(m = 2{,}000\), \(n = 5{,}120\), and \(b := A x_{\text{true}} + \varepsilon\), where \(\varepsilon \sim \mathcal{N}(0, 0.01)\).
The true signal \(x_{\text{true}}\) is a vector of zeros, except for \(100\) of its components.
We set \(\lambda = 0.1 \|A^T b \|_{\infty}\).
All algorithms start from the same randomly generated, hence non-sparse, \(x_0\).
For this problem, we compare R2 with R2DH variants (Spec, PSB, Andrei, and DBFGS).
The objective is to find a good subsolver for R2N, as the subproblem in Line~\ref{alg:R2N:step-computation} of \Cref{alg:R2N} has a structure similar to~\eqref{eq:bpdn}.

\begin{figure}[ht]%
  \centering
  \resizebox{.49\textwidth}{!}{%
    \includetikzgraphics{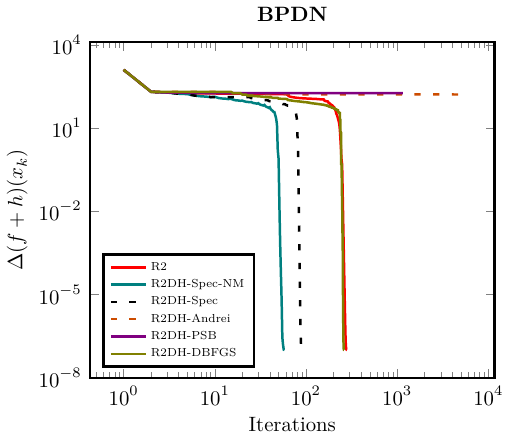}
  }
  \hfill
  \resizebox{.49\textwidth}{!}{%
    \includetikzgraphics{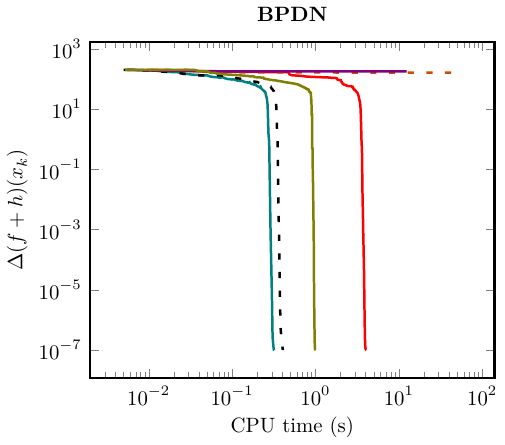}
  }
  \caption{BPDN objective vs.\ iterations (left) and CPU time (right).}%
  \label{fig:BPDN}
\end{figure}

\begin{table}[ht]%
  \footnotesize
  \centering
  \caption{Comparison of different solvers on the BPDN problem.}%
  \label{tab:bpdn}
  \begin{tabular}{rrrrrrrrr}
  \hline
  Solver & $f$ & $h/\lambda$ & $\Delta(f+h)$ & $\nu^{-1}\|s\|$ & $\#f$ & $\#\nabla f$ & \#prox & $t$($s$) \\\hline
  R2 & $9.86$e$-02$ & 100 & $4.26$e$-10$ & $1.6$e$-05$ & 281 & 273 & 280 & 4.08 \\
  R2DH-Spec-NM & $9.86$e$-02$ & 100 & $0.00$e$+00$ & $3.8$e$-06$ & 58 & 58 & 57 & 0.32 \\
  R2DH-Spec & $9.86$e$-02$ & 100 & $4.36$e$-10$ & $1.8$e$-05$ & 89 & 59 & 88 & 0.41 \\
  R2DH-Andrei & $3.00$e$+00$ & 2922 & $1.61$e$+02$ & $2.2$e$+00$ & 5001 & 4987 & 9991 & 55.49 \\
  R2DH-PSB & $1.34$e$-09$ & 3300 & $1.79$e$+02$ & $2.0$e$-05$ & 1153 & 1153 & 2311 & 12.36 \\
  R2DH-DBFGS & $9.86$e$-02$ & 100 & $1.24$e$-10$ & $1.3$e$-05$ & 262 & 153 & 261 & 1.00 \\\hline
\end{tabular}

\end{table}

\Cref{fig:BPDN} shows that all solvers reach similar accuracy, except for the PSB and Andrei variants.
R2DH-Spec-NM displays the best performance, followed by the R2DH-Spec and closely by R2DH-DBFGS variants, although it requires more evaluations to achieve stationarity.
\Cref{tab:bpdn} shows that all R2DH variants surpass R2 in all measures, except for R2DH-PSB and R2DH-Andrei, which either require more evaluations to attain the same level of accuracy or appear to converge to a different stationary point.
R2 and all other R2DH variants identify a similarly-sparse solution.
R2DH-Andrei requires significantly more evaluations and time than other R2DH variants and hits the iteration limit before~\eqref{eq:r2n-stop-test} is triggered.
Note that R2DH-DBFGS requires fewer function and gradient evaluations, as well as less time, than R2, but struggles to compete with R2DH-Spec-NM and R2DH-Spec.
R2DH-Spec-NM is more efficient than R2DH-Spec, it avoids the unsuccessful iterations that R2DH-Spec falls into.
Given its strong performance, in the following experiments we adopt R2DH-Spec-NM both as our default implementation of R2DH and as the R2DH subsolver.

\subsection{Matrix completion}

We address the matrix completion problem from~\citep{yu-zhang-2022} with rank nuclear norm regularizers to recover a low-rank matrix from noisy observations.
The problem is formulated as
\begin{equation}%
  \label{eq:mc-problem}
  \minimize{X} \ \tfrac{1}{2} \|P_{\Omega}
  (X - M)\|_F^2 + \lambda h(X),
\end{equation}
where \(X \in \mathbb{R}^{n \times n}\) and \(n=120\).
Here, \(\lambda=10^{-1}\) is a weight, and \(h(X)\) is either \(\rank(X)\) or \(\|X\|_*\), \(M\) is formed by applying a standard two-component Gaussian mixture model (GMM) to a low-rank matrix \(X_r\).
Specifically, \(M\) is computed as:
\[
  M = (1-c)(X_r + \mathcal{N}(0,\sigma_A^2)) + c(X_r + \mathcal{N}(0,\sigma_B^2)),
\]
where \(\mathcal{N}(0,\sigma_A^2)\) represents the noise component with variance \(\sigma_A^2\), and \(\mathcal{N}(0,\sigma_B^2)\) represents the influence of outliers with a larger variance \(\sigma_B^2\).
The parameter \(c\) controls the relative proportion of noise and outliers in the observed matrix \(M\).
Finally, \(P_\Omega\) is a linear operator that extracts entries \((i, j) \in \Omega\) and sets unobserved entries to zero, where
\[
  \Omega = \left\{ (i, j) \in \{1, \dots, m\} \times \{1, \dots, n\} \;\middle|\; R_{ij} < s_r \right\},
\]
\(R_{ij} \sim \mathcal{U}(0, 1)\), and \(s_r\) is a threshold that determines the sparsity of the observed matrix.

\begin{figure}[ht]%
  \resizebox{.49\textwidth}{!}{%
    \includetikzgraphics{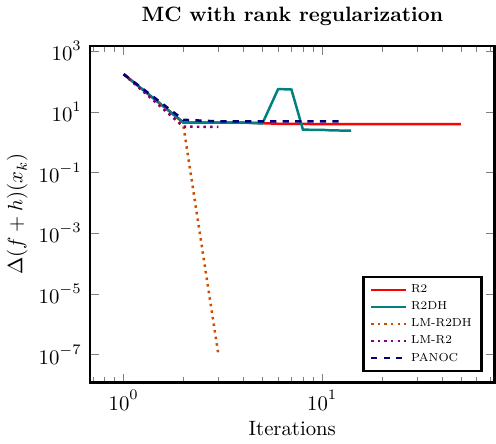}
  }
  \hfill
  \resizebox{.49\textwidth}{!}{%
    \includetikzgraphics{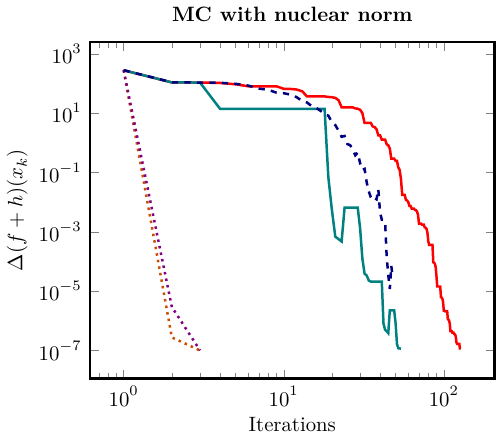}
  }
  \caption{Objectives vs.\ iterations for MC with rank (left) and nuclear norm (right) regularizers.}%
  \label{fig:comp}
\end{figure}

For all solvers, we select a random initial matrix and set the rank of \(X_r\) to \(40\).
The parameters are set to \(c = 0.2\), \(\sigma_A^2 = 0.0001\), \(\sigma_B^2 = 0.01\) and \(s_r = 0.8\).
Given that the smooth part of~\eqref{eq:mc-problem} is a linear least-squares residual, we apply the Levenberg-Marquardt (LM) algorithm from \citet[Algorithm~\(4.1\)]{aravkin-baraldi-orban-2024}, which is a specific instance of R2N with \(B_k = J_k^T J_k\), where \(J_k\) is the Jacobian of the least-squares residual at iteration \(k\).
Notably, R2DH can serve as a subproblem solver within LM---this combination is referred to as LM-R2DH, in contrast to the default LM-R2.
We compare the performance of R2, R2DH, LM-R2, LM-R2DH and PANOC in \Cref{fig:comp} and \Cref{tab:comp,tab:comp-nn}.
In \Cref{tab:comp,tab:comp-nn}, the column $\#\nabla f$ is replaced by the number of Jacobian or adjoint products $\#J$.

\begin{table}[ht]%
  \footnotesize
  \centering
  \caption{Comparison of different solvers for matrix completion problem with rank regularizer.}%
  \label{tab:comp}
  \begin{tabular}{rrrrrrrrr}
  \hline
  Solver & $f$ & $h/\lambda$ & $\Delta(f+h)$ & $\nu^{-1}\|s\|$ & $\#f$ & $\#J$ & \#prox & $t$($s$) \\\hline
  R2 & $1.61$e$-09$ & 111 & $3.90$e$+00$ & $1.7$e$-05$ & 50 & 43 & 49 & 0.24 \\
  R2DH & $1.09$e$-13$ & 96 & $2.40$e$+00$ & $4.7$e$-07$ & 14 & 14 & 13 & 0.07 \\
  LM-R2DH & $1.48$e$-11$ & 72 & $0.00$e$+00$ & $5.4$e$-06$ & 3 & 283 & 103 & 0.59 \\
  LM-R2 & $2.67$e$-11$ & 104 & $3.20$e$+00$ & $1.6$e$-06$ & 3 & 585 & 201 & 1.18 \\
  PANOC & $1.07$e$-09$ & 120 & $4.80$e$+00$ & $5.1$e$-06$ & 32 & 32 & 19 & 0.10 \\\hline
\end{tabular}

\end{table}

\begin{table}[ht]%
  \footnotesize
  \centering
  \caption{Comparison of different solvers for matrix completion problem with nuclear norm regularizer.}%
  \label{tab:comp-nn}
  \begin{tabular}{rrrrrrrrr}
  \hline
  Solver & $f$ & $h/\lambda$ & $\Delta(f+h)$ & $\nu^{-1}\|s\|$ & $\#f$ & $\#J$ & \#prox & $t$($s$) \\\hline
  R2 & $1.00$e$-02$ & $7.5$e$+00$ & $1.54$e$-08$ & $1.9$e$-05$ & 128 & 79 & 127 & 0.65 \\
  R2DH & $1.00$e$-02$ & $7.5$e$+00$ & $1.50$e$-08$ & $7.0$e$-06$ & 54 & 26 & 53 & 0.23 \\
  LM-R2DH & $1.00$e$-02$ & $7.5$e$+00$ & $0.00$e$+00$ & $3.2$e$-08$ & 3 & 589 & 201 & 1.08 \\
  LM-R2 & $1.00$e$-02$ & $7.5$e$+00$ & $1.07$e$-13$ & $5.9$e$-08$ & 3 & 555 & 201 & 1.05 \\
  PANOC & $1.00$e$-02$ & $7.5$e$+00$ & $2.83$e$-05$ & $8.0$e$-06$ & 103 & 103 & 55 & 0.41 \\\hline
\end{tabular}

\end{table}

\Cref{fig:comp} shows that LM-R2DH stands out in terms of final objective value for the rank regularizer followed by R2DH, while for the nuclear norm regularizer, all solvers achieve similar final objective value.
Variants of LM require the fewest objective evaluations, although they demand many Jacobian-vector products, as seen in \Cref{tab:comp,tab:comp-nn}.
For the rank regularizer, \Cref{tab:comp} shows that, while R2DH requires more objective evaluations than either LM variant, it performs significantly fewer Jacobian-vector products.
It is followed by PANOC, which requires more objective evaluations and Jacobian-vector products than R2DH.
In terms of the objective value, LM-R2DH outperforms all the other solvers, but at the cost of additional Jacobian-vector products and proximal evaluations, and provides the solution with the best objective value, and, in particular, the lowest-rank solution.
Finally, for the nuclear norm regularizer, LM variants behave almost identically according to \Cref{fig:comp} and \Cref{tab:comp-nn}.
R2DH outperforms R2 and PANOC in all measures, as shown in \Cref{tab:comp-nn}.

\subsection{General regularized problems}

In this section, we illustrate the performance of R2N on two test problems.
The first problem addresses an image recognition task using a support vector machine (SVM) similar to those in \citep{aravkin-baraldi-orban-2022}.
The objective is to use this nonlinear SVM to classify digits from the MNIST dataset, specifically distinguishing between ``1'' and ``7'', while excluding all other digits.
A sparse support is imposed using an \(\ell_0\) regularizer.
The optimization problem is given by
\[
  \minimize{x \in \R^n} \ \tfrac{1}{2} \|\mathbf{1} - \tanh(b \odot \langle A, x \rangle)\|^2 + \lambda \|x\|_0,
\]
where \(\lambda = 10^{-1}\) and \( A \in \mathbb{R}^{m \times n} \), with \( n = 784 \) representing the vectorized size of each image.
In our tests, we use the training dataset, which includes \( m = 13{,}007 \) images.
Here, \(\odot\) denotes the elementwise product between vectors, and \(\mathbf{1} = (1, \ldots, 1)\).

The second problem is from \citep{chouzenoux-martin-pesquet-2023,stella-themelis-patrinos-2017} and arises in image denoising and deblurring applications.
The related optimization problem is given by
\[
  \minimize{x \in \R^n} \ \sum_{i=1}^n \log \left( (A x - b)_i^2 + 1 \right) + \lambda \|x\|_1,
\]
where \( \lambda = 10^{-4} \) and \( A \in \R^{n \times n} \) with \( n = 256^2 \) is a Gaussian blur operator.
The term \( b \) denotes the blurred image with added Gaussian noise.
In our test, \( b \) is the blurred version of the cameraman image \( x^* \) with added Gaussian noise, i.e., \( b = A x^* + \text{noise} \).
The smooth part related to the two optimization problems is neither quadratic nor linear least squares, but a general non-convex problem.

We compare the performance of five methods: R2, R2DH, R2N-R2 (R2N with R2 as a subsolver), R2N-R2DH (R2N with R2DH as a subsolver) and PANOC \citep{stella-themelis-sopasakis-patrinos-2017}.

\begin{figure}[ht]%
  \resizebox{.48\textwidth}{!}{%
    \includetikzgraphics{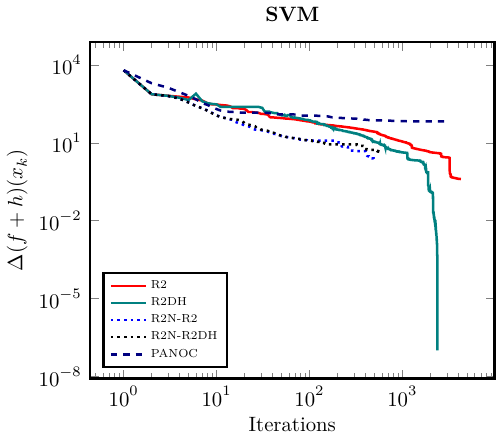}
  }
  \hfill
  \resizebox{.48\textwidth}{!}{%
    \includetikzgraphics{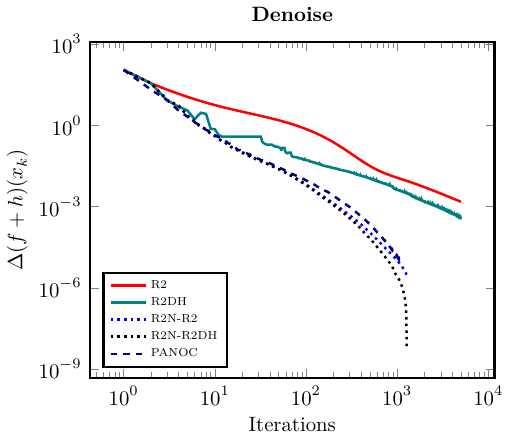}
  }
  \caption{Plots of the objective vs.\ iterations related to SVM (left) and denoise (right).}%
  \label{fig:general_pb}
\end{figure}

\begin{table}[h!]%
  \footnotesize
  \centering
  \caption{Comparison of different solvers on the nonlinear SVM problem.}%
  \label{tab:svm-l0}
  \begin{tabular}{rrrrrrrrr}
  \hline
  Solver & $f$ & $h/\lambda$ & $\Delta(f+h)$ & $\nu^{-1}\|s\|$ & $\#f$ & $\#\nabla f$ & \#prox & $t$($s$) \\\hline
  R2 & $1.40$e$+01$ & 122 & $4.02$e$-01$ & $2.0$e$-05$ & 4267 & 3303 & 4266 & 31.72 \\
  R2DH & $1.40$e$+01$ & 118 & $0.00$e$+00$ & $2.0$e$-05$ & 2384 & 1369 & 2383 & 15.12 \\
  R2N-R2 & $1.20$e$+01$ & 161 & $2.29$e$+00$ & $1.9$e$-05$ & 513 & 313 & 51201 & 4.29 \\
  R2N-R2DH & $1.60$e$+01$ & 144 & $4.59$e$+00$ & $2.0$e$-05$ & 561 & 297 & 54975 & 4.49 \\
  PANOC & $1.40$e$+01$ & 784 & $6.66$e$+01$ & $1.7$e$-05$ & 7338 & 7338 & 4211 & 73.84 \\\hline
\end{tabular}

\end{table}

\begin{table}[h!]%
  \footnotesize
  \centering
  \caption{Comparison of different solvers on the denoising problem.}%
  \label{tab:den}
  \begin{tabular}{rrrrrrrrr}
  \hline
  Solver & $f$ & $h/\lambda$ & $\Delta(f+h)$ & $\nu^{-1}\|s\|$ & $\#f$ & $\#\nabla f$ & \#prox & $t$($s$) \\\hline
  R2 & $5.89$e$-02$ & $3.6$e$+03$ & $1.51$e$-03$ & $9.8$e$-04$ & 5001 & 4998 & 5000 & 66.86 \\
  R2DH & $5.88$e$-02$ & $3.6$e$+03$ & $3.74$e$-04$ & $5.5$e$-04$ & 5001 & 2924 & 5000 & 51.25 \\
  R2N-R2 & $5.88$e$-02$ & $3.6$e$+03$ & $2.53$e$-06$ & $1.9$e$-05$ & 1327 & 1327 & 132601 & 325.60 \\
  R2N-R2DH & $5.88$e$-02$ & $3.6$e$+03$ & $0.00$e$+00$ & $2.0$e$-05$ & 1269 & 1269 & 126139 & 301.78 \\
  PANOC & $5.87$e$-02$ & $3.6$e$+03$ & $8.95$e$-06$ & $1.9$e$-05$ & 2192 & 2192 & 1104 & 36.39 \\\hline
\end{tabular}

\end{table}

As shown in \Cref{tab:den,tab:svm-l0}, for both problems, the R2N variants outperform R2, R2DH and PANOC in terms of objective and gradient evaluations, though they require more proximal operator evaluations.
Both R2N-R2DH and R2N-R2 have comparable performance and reach very good solutions compared to the other methods.

Note that for the non-linear SVM problem, as indicated in \Cref{fig:general_pb} and \Cref{tab:svm-l0}, although R2DH reduces the objective function the most, it requires a higher number of evaluations of $f$ and $\nabla f$ than both R2N-R2DH and R2N-R2.
PANOC is the least efficient method for this problem, as it requires the most evaluations and time while achieving the worst objective value.

For the denoising problem, the R2N variants require the fewest objective and gradient evaluations, but incur significantly more proximal-operator evaluations and longer runtime.
These costs arise primarily from solving the subproblem in Line~\ref{alg:R2N:step-computation} of \Cref{alg:R2N}.
Moreover, since the dimension of the denoising problem is \(65,536\), the cost of evaluating the proximal operator is not negligible.
PANOC requires roughly twice as many function and gradient evaluations as the R2N variants, while using fewer proximal-operator evaluations and less time.

\subsection{FitzHugh-Nagumo inverse problem}
We consider an inverse problem for recovering the parameters of a nonlinear ordinary differential equation (ODE) model.
Let \( x \in \mathbb{R}^{p} \) with \( p = 5 \) denote the model parameters, and let
\( F : \mathbb{R}^{p} \to \mathbb{R}^{2(N+1)} \) map \( x \) to the time series of state variables obtained by solving the FitzHugh--Nagumo (FH) neuron activation model \citep{aravkin-baraldi-orban-2022,aravkin-baraldi-orban-2024}
\begin{equation}
  \label{eq:fh-odes}
  \frac{\mathrm{d}
    V}{\mathrm{d}t} = \bigl(V - \tfrac{1}{3}V^{3} - W + x_1\bigr)\,x_2^{-1},
  \qquad
  \frac{\mathrm{d}W}{\mathrm{d}t} = x_2\,(x_3 V - x_4 W + x_5),
\end{equation}
initialized at \((V,W)=(2,0)\) and integrated over \(t \in [0,20]\) over a uniform grid of \(N\) steps.
We take \(N=1000\), i.e., \(\Delta t=0.02\)s.
We denote \(V(t_i;x)\) and \(W(t_i;x)\) the numerical solution at grid points \(\{t_i\}_{i=0}^{N}\), and set \(F(x) = \bigl(v(x),w(x)\bigr)\) with \(v(x)=(V(t_0;x),\ldots,V(t_N;x))\) and \(w(x)=(W(t_0;x),\ldots,W(t_N;x))\).
Synthetic observations are generated as \(b = F(x_{\textup{true}}) + \mathcal{N}(0,\,0.1^2 I)\), using the sparse ground truth \(x_{\textup{true}}=(0,1,0,0,0)\).
We estimate \(x\) using the sparsity-regularized problem
\begin{equation}%
  \label{eq:fh}
  \minimize{x\in\R^{p}} \ \tfrac{1}{2}\,\|F(x)-b\|_2^2 \;+\; \|x\|_0,
\end{equation}
where \(F\) and its Jacobian are computed by solving and performing automatic differentiation on~\eqref{eq:fh-odes}.

\begin{figure}[ht]%
  \centering
  \resizebox{.48\textwidth}{!}{%
    \includetikzgraphics{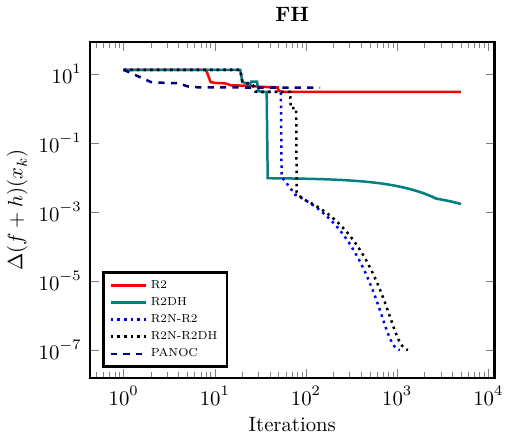}
  }
  \caption{FH objective vs.\ iterations.}%
  \label{fig:fh}
\end{figure}

\begin{table}[ht]%
  \footnotesize
  \centering
  \caption{Comparison of different solvers on the FH problem.}%
  \label{tab:fh}
  \begin{tabular}{rrrrrrrrr}
  \hline
  Solver & $f$ & $h/\lambda$ & $\Delta(f+h)$ & $\nu^{-1}\|s\|$ & $\#f$ & $\#\nabla f$ & \#prox & $t$($s$) \\\hline
  R2 & $1.17$e$+00$ & 4 & $3.00$e$+00$ & $3.5$e$-02$ & 5001 & 3758 & 5000 & 5.67 \\
  R2DH & $1.17$e$+00$ & 1 & $1.69$e$-03$ & $1.2$e$-02$ & 5001 & 2494 & 5000 & 7.20 \\
  R2N-R2 & $1.17$e$+00$ & 1 & $0.00$e$+00$ & $2.0$e$-05$ & 1067 & 1036 & 106601 & 1.38 \\
  R2N-R2DH & $1.17$e$+00$ & 1 & $7.26$e$-11$ & $2.0$e$-05$ & 1313 & 1270 & 125856 & 1.60 \\
  PANOC & $1.16$e$+00$ & 5 & $4.00$e$+00$ & $1.9$e$-05$ & 380 & 380 & 237 & 46.07 \\\hline
\end{tabular}

\end{table}

The FH problem is challenging due in particular to the fact that the gradient of the smooth term in~\eqref{eq:fh} is not Lipschitz continuous.

For this problem, the R2N variants and R2DH recover a sparse solution, as shown in \Cref{tab:fh,tab:fh-values}.
PANOC requires fewer objective, gradient, and proximal evaluations than R2N and R2DH, but struggles to converge to a sparse solution, is very slow, and requires the most time.
R2N-R2 (closely followed by R2N-R2DH) shows the best performance and successfully recovers a sparse solution.
Although the R2N variants require more proximal operator evaluations than the other methods, the cost of evaluating the proximal operator is negligible, since the FH problem is low-dimensional with only \(p=5\) parameters.

\begin{table}[ht]%
  \footnotesize
  \centering
  \caption{Solution recovered by the different solvers on the FH problem.}%
  \label{tab:fh-values}
  \begin{tabular}{rrrrrr}
  \hline
  True   & R2      & R2DH    & R2N-R2   & R2N-R2DH    & PANOC           \\\hline
  $0.00$ & $0.00$  & $0.00$  & $0.00$   & $0.00$      & $2.51$e$-10$    \\
  $1.00$ & $1.30$  & $1.47$  & $1.78$   & $1.78$      & $1.71$          \\
  $0.00$ & $-0.25$ & $0.00$  & $0.00$   & $0.00$      & $-263.31$       \\
  $0.00$ & $0.81$  & $0.00$  & $0.00$   & $0.00$      & $2181.53$       \\
  $0.00$ & $0.43$  & $0.00$  & $0.00$   & $0.00$      & $463.48$        \\\hline
\end{tabular}
\end{table}

\section{Discussion}%
\label{sec:discussion}

We proposed method R2N, a modified quasi-Newton method for nonsmooth regularized problems.
R2N generalizes both R2 \citep{aravkin-baraldi-orban-2022} and LM \citep{aravkin-baraldi-orban-2024} and enjoys convergence properties without assuming Lipschitz continuity of \(\nabla f\) or boundedness of the model Hessians.
Inspired by \citet{diouane-habiboullah-orban-2024}, who work on trust-region methods for smooth optimization, we propose a complexity analysis of R2N to handle potentially unbounded model Hessians.
Unlike traditional complexity analyses that assume uniformly bounded model Hessians, our study covers practical cases, including quasi-Newton updates such as PSB, BFGS, and SR1 by bounding the model Hessian growth with a power of the number of successful iterations---a reasonable bound as, in practice, it is uncommon to update quasi-Newton approximations on unsuccessful iterations.
Nevertheless, \citet{diouane-habiboullah-orban-2024} show that similar complexity bounds continue to hold when the model Hessians are bounded by a power of the number of iterations, and not just the number of successful iterations.
Because their analysis uses similar arguments, their complexity bounds continue to hold for R2N.

Numerical illustrations show the strong potential of our implementation of R2N and some of its variants (both as a main solver and as a subproblem solver) compared to PANOC \citep{stella-themelis-sopasakis-patrinos-2017}.
In particular, diagonal variants are competitive with, and often outperform, R2 when used as a subsolver inside R2N.
One of the main advantages of R2N in practice is that proximal operators are easier to compute than in TR \citep{aravkin-baraldi-orban-2022}.
We illustrated that advantage by solving rank and nuclear norm-regularized problems.
One way to further enhance the performance of R2N is to use more efficient subproblem solvers, for example, by generalizing those proposed in \citep{becker-fadili-2012,becker-fadili-ochs-2019,kanzow-lechner-2024} for convex \(h\) and reducing the required number of proximal operator evaluations.
Alternatively, in certain cases, the subproblems can be solved exactly, as in~\citep{diouane-gollier-orban-2024}.

R2N convergence analysis arguments can be used to update and strengthen the existing convergence analysis of methods R2, TR, TRDH \citep{leconte-orban-2025}, LM and LMTR \citep{aravkin-baraldi-orban-2024}.
In follow-up research, we aim to identify the nature of limit points under our assumptions.

\section*{Acknowledgement}
The authors would like to thank Maxence Gollier for his help with an improved implementation of R2N, as well as two anonymous reviewers for their constructive comments and suggestions that helped improve the paper.



\small
\bibliographystyle{abbrvnat}
\bibliography{abbrv,r2n-siam}
\normalsize

\appendix


\end{document}